\documentclass[10pt]{amsart}

\usepackage{comment}
\usepackage[latin2]{inputenc}
\usepackage{amsmath}
\usepackage{graphicx}
\usepackage{amssymb}
\usepackage{esint}
\usepackage{tikz}
\usepackage{xxcolor}
\usepackage{color}
\usepackage{amsthm}
\usepackage{epsfig}
\usepackage[english]{babel}
\usepackage{hyperref}
\hypersetup{
    colorlinks=true,
    linkcolor=blue,
    citecolor = red,
    pdfpagemode=FullScreen,
    }

\usepackage{lipsum}                     
\usepackage{xargs}                      
\usepackage[pdftex,dvipsnames]{xcolor} 

\usepackage[colorinlistoftodos,prependcaption,textsize=footnotesize]{todonotes}
\newcommandx{\fc}[2][1=]{\todo[linecolor=red,backgroundcolor=red!25,bordercolor=red,#1]{FC: #2}}
\newcommandx{\fcinline}[2][1=]{\todo[inline,linecolor=red,backgroundcolor=red!25,bordercolor=red,#1]{FC: #2}}

\newcommand{\tilt}[3]{\mathtt{TiltEx}_{{#2},{#3}}(#1)}
\newcommand{\fulltilt}[4]{\mathtt{TiltEx}[#4]_{{#2},{#3}}(#1)}
\newcommand{\hfterms}[3]{\int_{#1}^{#2}\mathtt{HFill}_{{#3}}}
\newcommand{\ratiodelta}{\Lambda}

\newtheorem{theorem}{Theorem}
\newtheorem{proposition}[theorem]{Proposition}
\newtheorem{lemma}[theorem]{Lemma}

\newtheorem{definition}[theorem]{Definition}
\newtheorem{remark}[theorem]{Remark}

\newtheorem*{theorem*}{Theorem}

\def\Xint#1{\mathchoice
{\XXint\displaystyle\textstyle{#1}}%
{\XXint\textstyle\scriptstyle{#1}}%
{\XXint\scriptstyle\scriptscriptstyle{#1}}%
{\XXint\scriptscriptstyle\scriptscriptstyle{#1}}%
\!\int}
\def\XXint#1#2#3{{\setbox0=\hbox{$#1{#2#3}{\int}$ }
\vcenter{\hbox{$#2#3$ }}\kern-.6\wd0}}

\def\dashint{\Xint-}

\allowdisplaybreaks[2]

\newcommand{\supp}{\operatorname{supp}}

\definecolor{Yellow}{rgb}{0.95,0.9,0.0} 
\definecolor{Red}{rgb}{0.8,0.1,0.1}
\definecolor{Green}{rgb}{0.1,0.65,0.2}
\definecolor{Blue}{rgb}{0.1,0.1,0.8}
\definecolor{Purple}{rgb}{0.7,0.1,0.7}
\definecolor{Grey}{rgb}{0.6,0.6,0.6}

\def\d{\, \mathrm{d}}
\def\R{\mathbb{R}}

\begin{document}

\newcommand{\InternalRemark}[1]{{\color{red}#1}}

\title[H\"older continuity for the 2d thin-film equation]{H\"older continuity of weak solutions to the thin-film equation in $d=2$}

\author{Federico Cornalba}
\address{University of Bath, Claverton Down, BA2 7AY, Bath, United Kingdom}
\email{fc402@bath.ac.uk}
\author{Julian Fischer}
\address{Institute of Science and Technology Austria, Am Campus 1, 3400 Klosterneuburg, Austria}
\email{julian.fischer@ista.ac.at}
\author{Erika Maringov\'a Kokavcov\'a}
\address{Institute of Science and Technology Austria, Am Campus 1, 3400 Klosterneuburg, Austria}
\email{erika.kokavcova@ist.ac.at}

\begin{abstract}
The thin-film equation $\partial_t u = -\nabla \cdot (u^n \nabla \Delta u)$ describes the evolution of the height $u=u(x,t)\geq 0$ of a viscous thin liquid film spreading on a flat solid surface. We prove H\"older continuity of energy-dissipating weak solutions to the thin-film equation in the physically most relevant case of two spatial dimensions $d=2$.
While an extensive existence theory of weak solutions to the thin-film equation was established more than two decades ago, even boundedness of weak solutions in $d=2$ has remained a major unsolved problem in the theory of the thin-film equation.
Due the fourth-order structure of the thin-film equation, De~Giorgi-Nash-Moser theory is not applicable.
Our proof is based on the hole-filling technique, the challenge being posed by the degenerate parabolicity of the fourth-order PDE.
\end{abstract}

\keywords{}

\maketitle

\section{Introduction}

The thin-film equation
\begin{align}
\label{tfe}
\partial_t u = -\nabla \cdot (u^n \nabla \Delta u)
\end{align}
(with $u:\mathbb{R}^d\times [0,T)\rightarrow \mathbb{R}_0^+$) describes the evolution of a viscous thin liquid film spreading on a solid surface. The parameter $n>0$ is related to the slip condition at the liquid-solid interface, with $n=3$ corresponding to a no-slip condition and with $n=2$ corresponding to a Navier slip condition.

The thin-film equation \eqref{tfe} may be regarded as the fourth-order analogue of the porous medium equation $\partial_t u=\Delta u^m=m\nabla \cdot (u^{m-1}\nabla u)$. While there are many similarities in the qualitative behavior of solutions -- such as preservation of nonnegativity of solutions, as well as the finite speed of propagation property of the free boundary $\partial \{x:u(x,t)>0\}$ -- there are important differences: Being a fourth-order equation, the thin-film equation lacks a comparison principle. Furthermore, as Remark~\ref{RemarkFailureHolderContinuity} below illustrates, for $n\geq 3$ one should not expect a regularizing effect of the evolution; in more than one space dimension, in this parameter regime solutions may never become H\"older continuous.

While the existence theory of weak solutions to the thin-film equation \eqref{tfe} is well-developed \cite{BerettaBertschDalPasso,BernisFriedman,Bernis,DalPassoGarckeGruen,
ElliottGarcke,ZAAGruen,Gruen2004} and the existence of strong solutions is known in many perturbative regimes around explicit solutions \cite{GiacomelliGnannKnuepferOtto,GiacomelliKnuepfer,GiacomelliKnuepferOtto,Gnann2}, regularity results for weak solutions have remained limited to those implied by the known integral estimates for the thin-film equation. Due to its fourth-order structure, De~Giorgi-Nash-Moser theory is not applicable to the thin-film equation. In the physically most relevant case of two spatial dimensions $d=2$, the energy estimate $\int |\nabla u(\cdot,t)|^2 \,\d x \leq \int |\nabla u_0|^2 \,\d x$ just barely fails to entail boundedness or H\"older continuity of weak solutions.  The question of boundedness or H\"older continuity of weak solutions in $d=2$ and for $n<3$ has remained one of the most important open problems in the theory of the thin-film equation.
 
In the present work, we prove the first regularity result for weak solutions to the thin-film equation beyond integral estimates: In two spatial dimensions $d=2$, we prove that any energy-dissipating weak solution $u$ as constructed in \cite{Gruen2004} is locally H\"older continuous in the sense $u\in C^\sigma_{loc}(\mathbb{R}^d\times (0,T))$ for some $\sigma=\sigma(n)>0$. If the initial data satisfy $\nabla u_0 \in L^p$ for some $p>2$, we even have global H\"older continuity up to the initial time $u\in C^\sigma(\mathbb{R}^d\times [0,T))$ for some $\sigma=\sigma(n,p)>0$.

\section{Overview of the Literature}

The mathematical theory of weak solutions to the thin-film equation relies on two basic integral estimates, the energy estimate
\begin{align}
\label{EnergyEstimate}
\partial_t \int \tfrac{1}{2} |\nabla u|^2 \d x = -\int u^n |\nabla \Delta u|^2 \d x
\end{align}
(which is formally readily verified by testing the PDE \eqref{tfe} with $\Delta u$)
and the so-called ``entropy estimates''
\begin{align}
\label{EntropyEstimate}
\partial_t \int \tfrac{1}{\alpha(1+\alpha)} u^{1+\alpha} \d x \leq -c(n,\alpha) \int |D^2 u^{(n+\alpha+1)/2}|^2 + |\nabla u^{(n+\alpha+1)/4}|^4 \d x
\end{align}
valid for any $\alpha\notin \{-1,0\}$ with $\tfrac{1}{2}<n+\alpha<2$. Unfortunately, no regularity results for weak solutions have been available that go beyond the regularity implied by \eqref{EnergyEstimate} and \eqref{EntropyEstimate}. As a consequence, also the uniqueness of weak solutions has remained a widely open problem.

This is despite the extensive existence theory for weak solutions to the thin-film equation having been developed several decades ago:
The first existence result for weak solutions was obtained by Bernis and Friedman \cite{BernisFriedman} in the case of one spatial dimension $d=1$; as their solution concept relies only on the energy estimate, it does not impose a contact angle condition at the free boundary and therefore suffers from a failure of uniqueness. In subsequent works by Beretta, Bertsch, and Dal~Passo \cite{BerettaBertschDalPasso} and Bertozzi and Pugh \cite{BertozziPughBlowup}, a notion of weak solution was developed for $n<2$ that is not subject to this immediate failure of uniqueness, based on the discovery of the family of entropy estimates \eqref{EntropyEstimate}; note that by the Morrey embedding, the regularity $(u^{(2-\delta+1)/4})_{x}\in L^4$ for a.\,e.\ $t>0$ entails a vanishing contact angle at the free boundary for a.\,e.\ $t>0$. Bernis \cite{Bernis} discovered the estimate $\int |(u^{(n+2)/6})_x|^6 + |(u^{(n+2)/3})_{xx}|^3 \d x \leq C \int u^n |u_{xxx}|^2 \d x$ valid for all smooth positive functions and for $n\in (\tfrac{1}{2},3)$, making it possible to develop a concept of weak solutions that relies on the energy estimate \eqref{EnergyEstimate} alone and that at the same time enforces a constraint of vanishing contact angle $|\nabla u|=0$ at the free boundary $\partial\{x:u(x,t)>0\}$.

In the case of multiple spatial dimensions, the first existence results for weak solutions were developed by Elliott and Garcke \cite{ElliottGarcke} and Gr\"un \cite{ZAAGruen}; for initial data with free boundary, these results were limited to $n<2$. Dal~Passo, Garcke, and Gr\"un \cite{DalPassoGarckeGruen} generalized the entropy estimate \eqref{EntropyEstimate} to the multidimensional case, allowing for the construction of weak solutions with contact angle constraint also for $n\in [2,3)$. Gr\"un \cite{GruenBernis} derived a multidimensional variant of the Bernis inequalities, enabling him to prove the existence of weak solutions subject to the energy dissipation property $\partial_t \int \tfrac{1}{2} |\nabla u|^2 \d x \leq -c\int |D^3 u^{(n+2)/2}|^2 + u^{n-2} |D^2 u|^2 |\nabla u|^2 + |\nabla u^{(n+2)/6}|^6 \d x$ in the parameter range $2-\sqrt{8/(8+d)}<n<3$ \cite{Gruen2004}.

For finite (nonzero) contact angle, only a limited number of existence results is available \cite{BertschGiacomelliKarali,Mellet,Otto}. The work by Otto \cite{Otto} exploits the Wasserstein gradient flow structure of the thin-film equation present in the particular case $n=1$; see also \cite{LisiniMatthesSavare} for a Wasserstein-like structure in the case of more general exponents $n<1$.

More recently, a rich theory of strong solutions in perturbative settings around explicitly known solutions has been developed, following the pioneering work by Giacomelli, Kn\"upfer, and Otto \cite{GiacomelliKnuepferOtto}, see 
\cite{GiacomelliGnannKnuepferOtto,GiacomelliKnuepfer,Gnann,Gnann2,
GnannIbrahimMasmoudi,GnannPetrache,John,Knuepfer,Knuepfer2,Seis} and the references therein.

Like the porous medium equation, the thin-film equation gives rise to a free boundary problem, the free boundary being the current boundary of the droplet $\partial \{x:u(x,t)>0\}$.
Unlike for the (second-order) porous medium equation, no comparison principle is available for the (fourth-order) thin-film equation; thus, the analysis of the qualitative behavior of solutions had to rely on localized versions of the energy and entropy estimates.
Finite speed of propagation of the free boundary -- in the sense that $\partial \{x:u(x,t)>0\}$ may only expand with a H\"older-like speed -- has been shown for the various settings of weak solutions \cite{BernisFinite,BernisFinite2,HulshofShishkov,BertschDalPassoGarckeGruen,
GruenOptimalRatePropagation}.
For $n>\tfrac{3}{2}$ it has been shown that the support of solutions cannot shrink \cite{BernisFriedman,BerettaBertschDalPasso,BertschDalPassoGarckeGruen}; in the regime $n<\tfrac{3}{2}$, the explicit solution to the thin-film equation $u(x,t):=(x+c_n t)_+^{3/n}$ demonstrates that the free boundary may recede, in stark contrast to the case of the porous medium equation.
In the particular case $d=1$, $n=1$, convergence to a self-similar solution was established by Carrillo and Toscani \cite{CarrilloToscani}.
In \cite{DalPassoGiacomelliGruen,GiacomelliGruen}, sufficient criteria for a waiting time phenomenon were established: If the initial data $u_0$ are flat enough near the free boundary $\partial\{x:u_0(x)\}>0$, the free boundary was shown to not move forward for some finite time before it could start advancing. Lower bounds on the propagation of the free boundary as well as upper bounds on waiting times have been established by the second author \cite{FischerARMA,FischerJDE,FischerAHP,DeNittiFischer}, based on the discovery of a new family of monotone quantities of the form $\int u^{1+\alpha} |x-x_0|^{-\gamma} \,\d x$ for suitable $\alpha\in (-1,0)$ and $\gamma>0$.

The thin-film equation has received considerable attention in physics; we only refer to \cite{Greenspan,OronDavisBankoff} regarding the classical thin-film equation, as well as to \cite{DavidovitchMoroStone,StochasticThinFilm} for its stochastic variant that incorporates thermal fluctuations in microscopic thin liquid films.
In recent years, the mathematical analysis of the latter has seen substantial developments; we refer to \cite{AgrestiSauerbrey,DareiotisGessGnannGruen,FischerGruen,
GessGvalaniKunickOtto,GvalaniTempelmayr,MetzgerGruen,sauerbrey2025solutions} and the references therein.

The hole-filling technique was originally developed in the context of second-order elliptic systems \cite{WidmanHoleFilling}; its applications in the parabolic context so far have been limited to strictly parabolic PDEs  \cite{StruweHoleFilling,FrehseSpecoviusNeugebauer}. In particular, the approach by Struwe \cite{StruweHoleFilling} appears inherently limited to strictly parabolic equations.



\section{Main Results}

In this work, we prove H\"older continuity of weak solutions to the thin-film equation \eqref{tfe} in the case of two spatial dimensions $d=2$, see Theorem~\ref{MainResult}. We will achieve this goal for the notion of weak solutions introduced by Gr\"un \cite{Gruen2004}, which are characterized by an energy dissipation inequality and are available in the parameter regime $2-\sqrt{4/5}\approx 1.106<n<3$.
We emphasize that this covers the majority of the parameter regime for which H\"older regularity of solutions may be expected: For parameter values $n\geq 3$, the so-called no-slip paradox is conjectured to prevent the motion of the contact line $\partial \{x:u(x,t)>0\}$; as shown in Remark~\ref{RemarkFailureHolderContinuity} below, this would immediately enable the construction of counterexamples to H\"older continuity. On the other hand, the lower bound on $n$ is expected to be technical and is a consequence of the current range of availability of the Bernis-Gr\"un inequalities (see Proposition~\ref{BernisGruenIneq}); an extension of the Bernis-Gr\"un inequalities to smaller values of $n$ would immediately imply a corresponding extension of our result.

\begin{theorem}[H\"older continuity of energy-dissipating weak solutions to the thin-film equation for $d=2$]
\label{MainResult}
Let $d=2$ and let $2-\sqrt{4/5}<n<3$. Let $u_0\in L^1(\mathbb{R}^2)\cap H^1(\mathbb{R}^2)$ be compactly supported. Let $u$ be a weak solution to the thin-film equation in the sense of Definition~\ref{DefinitionEnergyDissipatingWeakSolution} below.
\begin{itemize}
\item[a)] Then $u$ is H\"older continuous for $t>0$ in the sense $u\in C^\sigma_{loc}(\mathbb{R}^2 \times (0,T))$ for some $\sigma=\sigma(n)>0$.
\item[b)] If the initial data have the additional regularity $\nabla u_0 \in L^p(\mathbb{R}^2)$ for some $p>2$, the solution $u$ is H\"older continuous up to the initial time in the sense $u\in C^\sigma(\mathbb{R}^2\times [0,T))$ with $\sigma=\sigma(n,p)>0$.
\end{itemize}
\end{theorem}

Our main result applies to the following notion of weak solution, which we shall call energy-dissipating weak solutions\footnote{In the original work \cite{Gruen2004} these solutions are referred to as ``strong solutions'' in order to distinguish it from the weaker notion of distributional solution; to avoid confusion with the more recent works on solutions to the thin-film equation in H\"older spaces \cite{GiacomelliKnuepferOtto,GiacomelliGnannKnuepferOtto,Gnann2}, we shall instead refer to these solutions as energy-dissipating weak solutions.}. Note that for any $u_0\in H^1(\mathbb{R}^d)$ with compact support and any $d$, $n$ in the parameter range stated below, existence of such energy-dissipating weak solutions has been shown by Gr\"un \cite{Gruen2004}.

\begin{definition}[Energy-dissipating weak solutions, see Definition~1.1 in \cite{Gruen2004}]
\label{DefinitionEnergyDissipatingWeakSolution}
	Let $d\in \{1,2,3\}$ and $n \in (2-\sqrt{8/(8+d)},3)$. Let $T>0$ and let $u_0\in \smash{H^1(\R^d)}$ have compact support.
	We call a nonnegative function $\smash{u\in L^\infty([0,T); H^1(\R^d)\cap L^1(\R^d))}$, $u\geq 0$, an \emph{energy-dissipating weak solution of the thin-film equation with zero contact angle and initial data $u_0$} if the following conditions are satisfied:
	\begin{itemize}
		\item[a)] We have $\nabla u^{\frac{n+2}{6}} \in L^6(\R^d \times [0,T))$, $u^{\frac{n-2}{2}}\nabla u \otimes D^2 u \in L^2(\R^d \times [0,T))$, and $\chi_{\{u > 0\}}u^{\frac{n}{2}}\nabla\Delta u \in L^2(\R^d \times [0,T))$.
		\item[b)] For all $\alpha \in (\max\left\{-1,\frac{1}{2} - n\right\}, 2-n) \setminus \{0\}$, we have $D^2 u^{\frac{1+n+\alpha}{2}} \in L^2(\R^d \times [0,T))$ and $\nabla u^{\frac{1+n+\alpha}{4}} \in L^4(\R^d \times [0,T))$.
		\item[c)] It holds that $u \in H^1([0,T);(W^{1,p}(\R^d))')$ for all $p > \frac{4d}{2d + n(2-d)}$.
		\item[d)] For any $\psi \in L^2([0,T),W^{1,\infty}(\R^d))$ and any $T>0$, we have
		\begin{equation}\label{def:endisWS}
		\int_{0}^{T} \langle \partial_t u, \psi \rangle_{(W^{1,p}(\R^d))'\times W^{1,p}(\R^d)} \d t = \int_0^T \int_{\R^d \cap \{u>0 \}} u^n \nabla \Delta u \cdot \nabla \psi \d x \d t.
		\end{equation}
		\item[e)] $u$ attains its initial data $u_0$ in the sense $\lim\limits_{t\to 0}u(\cdot,t)=u_0(\cdot)$ in $L^1(\R^d)$.
	\end{itemize}
\end{definition}

Note that the upper bound $n<3$ for our H\"older continuity result is expected to be optimal:
\begin{remark}
\label{RemarkFailureHolderContinuity}
For $n\geq 3$, it is conjectured that the support of (suitably defined) weak solutions to the thin-film equation should remain constant over time. Assuming this conjecture, we could immediately construct a counterexample to H\"older continuity of solutions: If the support of solutions remains constant in time, this allows us to obtain a new solution by adding any two solutions with strictly disjoint initial support. Choose any nonnegative $\eta \in C^\infty_{cpt}(B_1 \setminus B_{1/2})$ on the annulus $B_1 \setminus B_{1/2}$, set $u_{0,m}(x):=m^{-1} \eta(4^{m} x)$, and consider weak solutions $u_m(x,t)$ with initial data $u_{0,m}(x)$. Then we may obtain a solution $u(x,t)$ with the initial data $u_0(x)=\sum_{m=1}^\infty u_{0,m}(x)$ as the sum $u(x,t)=\sum_{m=1}^\infty u_m(x,t)$; note that due to $\sum_{m=1}^\infty m^{-2}<\infty$ we have $u_0\in H^1(\mathbb{R}^2)$. By preservation of mass, for any $t>0$ and any $m$ we may then find $x\in B_{4^{-m}}$ with $u_m(x,t)\geq c m^{-1}$, which entails $|u(x,t)-u(0,t)|=|u(x,t)|\geq cm^{-1}\geq c |\log |x||^{-1}$.
\end{remark}
At the level of our method, our proof strategy indeed breaks down for $n\geq 3$: It is known that for $n\geq 3$ no Bernis-type estimate of the form $\int |\nabla u^{(n+2)/6}|^6 \d x \leq C \int u^n |\nabla \Delta u|^2 \d x$ can hold, thereby naturally limiting our approach to $n<3$.

\section{Outline of the strategy}

\subsection{The hole-filling technique: The elliptic setting}
Classically, the hole-filling technique provides a quick proof of H\"older continuity of weak solutions to linear elliptic systems in $d=2$ \cite{WidmanHoleFilling}.
Consider a weak solution $u$ to the linear elliptic system $-\nabla \cdot (a\nabla u)=0$, where $a(x)$ is uniformly elliptic and bounded. Choose a smooth cutoff $\eta\geq 0$ with $\eta\equiv 1$ in $B_r$ and $\eta\equiv 0$ outside of $B_{2r}$ as well as with $|\nabla \eta|\leq r^{-1}$. Testing with $(u-b)\eta^2$ for any constant $b$ yields the energy estimate (Caccioppoli inequality)
\begin{align*}
\int_{B_r} |\nabla u|^2 \,\d x
\leq C r^{-2} \inf_{b\in \mathbb{R}}  \int_{B_{2r}\setminus B_r} |u-b|^2 \,\d x.
\end{align*}
By the Poincar\'e inequality, this simplifies to $\int_{B_r} |\nabla u|^2 \,\d x
\leq \tilde C \int_{B_{2r}\setminus B_r} |\nabla u|^2 \,\d x$. Adding $\tilde C \int_{B_r} |\nabla u|^2 \,\d x$ to both sides (``filling the hole'' in the integral on the right-hand side) yields the bound
\begin{align*}
(1+\tilde C) \int_{B_r} |\nabla u|^2 \,\d x
\leq \tilde C \int_{B_{2r}} |\nabla u|^2 \,\d x.
\end{align*}
An iteration of this estimate then yields the existence of a small exponent $\sigma>0$ with
\begin{align*}
\int_{B_r} |\nabla u|^2 \,\d x
\leq C \bigg(\frac{r}{R}\bigg)^{2\sigma} \int_{B_R} |\nabla u|^2 \,\d x
\quad\text{for any }R>0\text{ and any }0<r<R,
\end{align*}
an estimate that is sufficient to establish H\"older continuity in $d=2$: The Poincar\'e inequality directly implies that $u$ belongs to the Campanato space $\mathcal{L}^{2+2\sigma,2}(\mathbb{R}^d)$, which in $d=2$ entails H\"older continuity. For the standard argument for this last step, we refer to the proof of Theorem~\ref{MainResult} below.

\subsection{Hole filling for uniformly parabolic equations}\label{subsec_unif_parab}
We next discuss the application of the hole-filling technique to uniformly parabolic PDEs with a structure similar to that of the thin-film equation, for now avoiding the critical issue of degenerate ellipticity in the thin-film equation \eqref{tfe}.
We emphasize that due to the specific fourth-order structure of the PDE, even in the uniformly parabolic setting our hole-filling approach differs from the more classical parabolic hole-filling approach by Struwe \cite{StruweHoleFilling}; in fact, it is closer in spirit to that of Frehse and Specovius-Neugebauer \cite{FrehseSpecoviusNeugebauer}, but involves testing with $\Delta u$ instead of $\partial_t u$.

We consider $\partial_t u = -\nabla \cdot (a(x) \nabla \Delta u)$, where $a(x)$ 
satisfies $1 \leq a(x) \leq a_{\text{max}}$. 
As will become apparent througout the paper, the tilt-excess-type quantity defined as
\begin{align}\label{abstract_def_tilt}
\fulltilt{t}{r}{\tilde{r}}{u} := \frac{1}{2}\int_{B_r} |\nabla u(x,t)-b_{\tilde{r}}(t) \cdot x -c_{\tilde{r}}(t)|^2 \,\d x,
\end{align}
where $b_{\tilde{r}}(t)$ (respectively, $c_{\tilde{r}}(t)$) is a suitable smoothed weighted average of second (respectively, first) derivatives of $u$ over $B_{2\tilde{r}}\setminus B_{\tilde{r}}$, will play a crucial role. 
\begin{remark}
As it will always be apparent which function $u$ we refer to, we will shorten $\fulltilt{t}{r}{\tilde{r}}{u}$ to $\tilt{t}{r}{\tilde{r}}$. Furthermore, for notational convenience, we may sometime omit time-dependencies.
\end{remark}

In order to estimate the tilt-excess-type quantity \eqref{abstract_def_tilt}, we consider a cutoff $\eta$ supported in $B_{2r}$ and equal to one in $B_r$, and perform the computation
\begin{align*}
& \frac{1}{2}\partial_t \int | \nabla u-b_r \cdot x -c_r |^2 \eta^6 \,\d x \\
& \quad  = \int \partial_t (\nabla u-b_r \cdot x -c_r) \cdot (\nabla u-b_r \cdot x -c_r)  \eta^6 \d x \\
& \quad = \int (\nabla (-\nabla \cdot (a\nabla\Delta u))-b_r' \cdot x -c_r') \cdot (\nabla u-b_r \cdot x -c_r)  \eta^6 \d x
\end{align*}
for some appropriate $b_r,c_r$ to be chosen.
Integrating by parts so as to remove all derivatives of order four and five, one gets
\begin{align*}
 & \frac{1}{2}\partial_t \!\! \int |\nabla u-b_r \cdot x -c_r|^2 \eta^6 \,\d x \\
 & \quad =
- \int a |\nabla \Delta u|^2 \eta^6\d x  
\\& \quad\quad\quad
- \int a\nabla \Delta u \cdot ((\Delta u - \operatorname{tr} b_r)\operatorname{Id} + D^2 u - b_r) 6\eta^5\nabla \eta \d x \\
& \quad \quad \quad 
- \int a\nabla \Delta u \cdot (\nabla u-b_r \cdot x -c_r)[30\eta^4|\nabla \eta|^2 + 6\eta^5\Delta \eta] \d x \\
& \quad \quad \quad - \int (b_r'\cdot x + c_r')\cdot (\nabla u - b_r\cdot x - c_r)\eta^6 \d x =: \sum_{i=1}^{4}{T_i}.
\end{align*}
We shift $T_1$ to the left-hand-side. Then, using the bound $a(x)\leq a_{\text{max}}$, we can bound $T_2$ and $T_3$ using Young's inequality, so as to detach and absorb the quantity $\int a|\nabla\Delta u|^2\eta^6\d x$ in the left-hand-side. This gives
\begin{align*}
&\partial_t \int |\nabla u-b_r \cdot x -c_r|^2 \eta^6 \,\d x
+\int a|\nabla \Delta u|^2 \eta^6 \,\d x
\\&
\leq C r^{-2} \int_{B_{2r}\setminus B_r} |D^2 u-b_r|^2 \,\d x + C r^{-4} \int_{B_{2r}\setminus B_r} |\nabla u-b_r\cdot x -c_r|^2 \,\d x
\\&~~~
\quad \quad \quad - \int (b_r'\cdot x + c_r')\cdot (\nabla u - b_r\cdot x + c_r)\eta^6 \d x.
\end{align*}
By properly tuning $b_r,c_r$, we may apply the Poincar\'e inequality to the first two terms on the right-hand side, thus obtaining
\begin{align}\label{pre_hole_fill_unif_parab}
&\partial_t \int |\nabla u-b_r \cdot x -c_r|^2 \eta^6 \,\d x
+\int_{B_r} |\nabla \Delta u|^2  \,\d x
\\&
\leq C \int_{B_{2r}\setminus B_r} |D^3 u|^2 \,\d x 
 - \int (b_r'\cdot x + c_r')\cdot (\nabla u - b_r\cdot x - c_r)\eta^6 \d x.\nonumber
\end{align}
where -- crucially -- we have used the non-degeneracy condition $1\leq a(x)$ in the left-hand-side above. 

In order to perform hole-filling in \eqref{pre_hole_fill_unif_parab}, we need a positive multiple of $\int_{B_r}{|D^3 u|^2\d x}$ on the left-hand side, and we may allow for a term of the form $\int_{B_{2r}}{|\nabla \Delta u|^2\d x}$ with sufficiently small pre-factor on the right-hand side. 
This can obtained by adding a multiple of the inequality $\int_{B_r}{|D^3 u|^2\d x} \lesssim \int_{B_{2r}}{|\nabla \Delta u|^2 \d x} + \int_{B_{2r}\setminus B_r}{|D^3 u|^2\d x}$ (see Lemma \ref{ThirdDerivativeBound}) with sufficiently small multiplicative factor.
\begin{remark}\label{rem_full_third_derivative}
Note that Lemma \ref{ThirdDerivativeBound} gives the means of controlling the full third derivative $D^3 u$ in hole-filling estimates: control of this quantity is not directly obtained from energy estimates (as opposed to control of $\nabla\Delta u$), but it nevertheless appears on the right-hand side of our estimates due to the use of Poincar\'e-type estimates.
\end{remark}

In addition, it can be shown (via calculations whose precise details we defer to later points in the paper) that the final term on the right-hand-side of \eqref{pre_hole_fill_unif_parab} can be dealt with without introducing any terms other than -- essentially -- those already present in the estimate. This yields
\begin{align*}
& \int_{B_{r}} \frac12|\nabla u-b_r \cdot x -c_r|^2(\cdot,t_2) \d x
+ \tilde c \int_{t_1}^{t_2} \!\!\! \int_{B_{r}} |\nabla \Delta u|^2 + |D^3 u|^2 \d x \d t
\\& \nonumber
\leq
\int_{B_{2r}} \frac12|\nabla u-b_r \cdot x -c_r|^2(\cdot,t_1) \d x
+\tilde C \int_{t_1}^{t_2} \!\!\! \int_{B_{2 r}\setminus B_r} \!\!\! |\nabla \Delta u|^2 + |D^3 u|^2 \d x \d t
\end{align*}
for some $\tilde c,\tilde C>0$.
Ultimately, one obtains the following hole-filling type estimate for the tilt-excess quantity \eqref{abstract_def_tilt}
\begin{align}\label{hole_filling_unif_parab}
&\int_{B_{r}} \frac12|\nabla u-b_r \cdot x -c_r|^2(\cdot,t_2) \d x
+ \hat C \int_{t_1}^{t_2} \!\!\! \int_{B_{r}} |\nabla \Delta u|^2 + |D^3 u|^2 \d x \d t
\\& \nonumber
\leq
\int_{B_{2r}} \frac12|\nabla u-b_r \cdot x -c_r|^2(\cdot,t_1) \d x
+(1-\theta)\cdot \hat C \int_{t_1}^{t_2} \!\!\! \int_{B_{2 r}} \!\!\! |\nabla \Delta u|^2 + |D^3 u|^2 \d x \d t
\end{align} 
with $\hat C:=\tilde c+\tilde C$ and with $1-\theta := \tfrac{\tilde C}{\tilde C + \tilde c} \in(0,1)$, or, with more succinct notation,
\begin{align*}
\tilt{t_2}{r}{r} + \hfterms{t_1}{t_2}{r} \leq \tilt{t_1}{2r}{r} + (1-\theta)\cdot  \hfterms{t_1}{t_2}{2r},
\end{align*}
where the term $\hfterms{t_1}{t_2}{r} = \hat C \int_{t_1}^{t_2} \!\!\! \int_{B_{r}} |\nabla \Delta u|^2 + |D^3 u|^2 \d x \d t$ contains all quantities which are involved in the hole-filling procedures above (in this uniformly parabolic example, $|\nabla \Delta u|^2$ and $|D^3 u|^2$).
This last inequality is the key to deducing $C^\sigma$ regularity in time and space for the solution $u$: Iterating this estimate, one deduces the excess-decay
\begin{align*}
\tilt{t_2}{r}{r} = \int_{B_{r}} \tfrac12|\nabla u-b_r \cdot x -c_r|^2(\cdot,t_2) \d x \lesssim r^\beta
\end{align*}
for some $\beta>0$ and thus spatial H\"older continuity. Getting an expression analogue to \eqref{hole_filling_unif_parab} in the case of the thin-film equation \eqref{tfe} is the core component of the paper.

\subsection{Hole-filling for the 2D thin-film equation} We now turn to the thin-film equation \eqref{tfe}. 
The key challenge (as compared to Subsection \ref{subsec_unif_parab}) is that the uniformly bounded term $a(x)$ is now replaced by the degenerate mobility term $u^n$, which -- a priori -- may come arbitrarily close to (or equal to) zero, or get arbitrarily large. Therefore, we need to adapt the estimates. Our crucial insight is to consider different hole-filling quantities for different times, depending on the categorisation according to the following definition.
\begin{definition}[Good and bad times]\label{DefnGoodBad}
Let $u$ be the solution to the thin-film equation \eqref{tfe} as per Definition \ref{DefinitionEnergyDissipatingWeakSolution}.

$\bullet$ A time $t\in[0,T]$ is said to be \textbf{good} over the ball $B_r$ (or, in short, for radius $r$) if the following `uniform-parabolicity'-type property holds:
\begin{equation}\label{goodt}
\sup_{x \in B_r} u(x,t) \leq 2 \inf_{x \in B_r} u(x,t).
\end{equation} 

$\bullet$ A time $t\in[0,T]$ is said to be \textbf{bad} over the ball $B_r$ (or, in short, for radius $r$) if the opposite holds, i.e., if
\begin{equation}\label{badt}
\inf_{x \in B_r} u(x,t) < \frac12 \sup_{x \in B_r} u(x,t).
\end{equation}
\end{definition}
\begin{remark}\label{remark_good_times}
Note that upon decreasing $r$ a good time $t$ always remains a good time, while a bad time $t$ may remain a bad time or may turn into a good time. 
\end{remark}

\emph{Key insights for analysis of good times}.
For good times, hole-filling estimates involve -- loosely speaking -- the same differential operators as in the uniformly parabolic case, simply with $u^n$ as a multiplier. This means that, for good times, hole-filling estimates involving the quantities 
$$
u^n |\nabla \Delta u|^2\mbox{ and }u^n |D^3 u|^2
$$ 
can be produced. 
In analogy to the parabolic case (see Remark \ref{rem_full_third_derivative} also) the term $u^n |\nabla \Delta u|^2$ naturally arises from basic manipulations of the thin-film equation, and can be dealt with relative ease. On the other hand, thanks to the good-time condition \eqref{goodt} (which -- very loosely speaking -- allows to treat the mobility as if it were `constant'), the term $u^n |D^3 u|^2$ can be included in the hole-filling estimates with an  estimate (see \eqref{second_formula}) which is -- in spirit --  analogous to the one of Lemma \ref{ThirdDerivativeBound}, which we have already mentioned. 
Furthermore, for good times, we can make use of Poincar\'e estimates involving the full derivatives $D^2 u$ and $D^3 u$ (see Lemma \ref{lemma_poincare}). This can be done thanks to the integrability properties recalled in Lemma \ref{finite_2_3_good_times}.

\emph{Key insights for analysis of bad times}.
For bad times, the third-derivative term $\int_{B_{2r}\setminus B_r} |D^3 u|^2 \d x$ that would arise from Poincar\'e-type estimates for the right-hand side of our energy estimate can not be controlled due to the possible degeneracy of the mobility $u^n$. Instead, we produce hole-filling estimates for the quantities 
$$
u^n |\nabla \Delta u|^2,\quad |\nabla u^{(n+2)/6}|^6.
$$
The distinctive quantity $\int_{B_{2r}} |\nabla u^{(n+2)/6}|^6 \d x$, as we will detail throughout the paper, stems from uses of the Bernis-Gr\"un inequality (see Lemma \ref{BernisGruenIneq}) and enables us to control the supremum $\sup_{B_{2r}} u^{n+2}$ for bad times via Morrey's inequality.

\emph{Combining estimates for good and bad times}. We obtain the hole-filling estimate \eqref{hole_filling_total} (generalising \eqref{hole_filling_unif_parab} from the uniformly parabolic case), and which can be succintly written as 
\begin{align}\label{hf_est_succinct_tfe}
& \tilt{t_2}{\delta r}{r} + \hfterms{t_1}{t_2}{\delta r}  \leq \tilt{t_1}{2r}{r} + (1-\theta) \cdot \hfterms{t_1}{t_2}{2r},
\end{align}
for some $\theta = \theta(n) \in (0,1)$, and some sufficiently small $\delta \in (0,1)$, where now
\begin{align}\label{hfill_def_tfe}
\hfterms{t_1}{t_2}{r} & = C_{\text{good}}\int_{t_1}^{t_2} \chi_{\text{good time}(r)}(t) \int_{B_{r}} u^n |\nabla \Delta u|^2 + u^n |D^3 u|^2 \d x \d t
\\&~~~ \nonumber
+ C_{\text{bad}}\int_{t_1}^{t_2} \chi_{\text{bad time}(r)}(t) \int_{B_{r}} u^n |\nabla \Delta u|^2 + |\nabla u^{(n+2)/6}|^{6} \d x \d t
\end{align}
for some positive constants $C_{\text{good}},C_{\text{bad}}$ only depending on $n$ and $d$.
The full details of the hole-filling estimate \eqref{hf_est_succinct_tfe} are given in Lemma \ref{LemmaHoldFill_TFE}. This key estimate then enables the derivation of a decay estimate for the tilt-excess $\tilt{t_2}{r}{r}\lesssim r^\beta$, thus establishing spatial H\"older continuity uniformly in time. An additional interpolation argument yields space-time H\"older continuity of the solution to the 2D thin-film equation (see our main result, Theorem \ref{MainResult}).

\section{A hole-filling estimate for the thin-film equation}


In this section, we consider the thin-film equation \eqref{tfe}, and prove the hole-filling estimate \eqref{hf_est_succinct_tfe}. In order to do this, we first need to rigorously define the quantities $b_r$ and $c_r$ that we have so far only colloquially introduced when defining the tilt-excess quantity $\fulltilt{t}{r}{\tilde{r}}{u}$, see \eqref{abstract_def_tilt}.


\begin{definition}[Smooth averaged second and first derivatives of thin-film solution over annuli]\label{Defn_b_c} Let $u$ be the solution to \eqref{tfe} in the sense of Definition \ref{DefinitionEnergyDissipatingWeakSolution}.
Let $\tilde{\eta}$ be a radially symmetric cutoff supported in $B_{2}\setminus B_1$ and equal to one in $B_{5/3}\setminus B_{4/3}$.
We define ``smoothly averaged second and first derivatives'' of the function $u$ in $B_{2r}\setminus B_r$, denoted by $b_r(t)$ (respectively, $c_r(t)$) as
\begin{align}
b_r(t)&:= -\frac{1}{\int_{B_{2r}} \tilde{\eta}(\frac{x}{r})\d x} \int_{B_{2r}} \nabla \Big(\tilde{\eta}\left(\frac{x}{r}\right) \Big) \otimes \nabla u \d x, \label{defb}\\ 
c_r(t)&:= \frac{1}{\int_{B_{2r}} \tilde{\eta}(\frac{x}{r})\d x} \int_{B_{2r}} \tilde{\eta}\left(\frac{x}{r}\right)\nabla u \d x. \label{defc}
\end{align}
\end{definition}
Note that $b_r(t)$ and $c_r(t)$ are well-defined (regardless or whether $t$ is a good or a bad time for $B_{2r}$) as 
$\smash{u\in L^\infty([0,T); H^1(\R^d)\cap L^1(\R^d))}$ (see Definition \ref{DefinitionEnergyDissipatingWeakSolution}). 
\begin{remark}
Besides being well-defined for any time $t$, if, in addition, $t$ is a good time for radius $2r$, the quantities $b_r(t)$ and $c_r(t)$ also allow for the use of Poincar\'e inequalities to bound quantities like $\nabla u(\cdot,t)-b_r(t)\cdot x-c_r(t)$ or $D^2 u(\cdot,t)-b_r(t)$, despite the fact that they do not coincide with the usual average as per standard Poincar\'e (see Lemma \ref{lemma_poincare} in Appendix \ref{appendix_Bs_Cs}).
We use $b_r$ and $c_r$ (as opposed to the standard spatial averages) due to their regularity properties granted by the kernel $\tilde{\eta}$ (in particular, lack of boundary terms when integrating by parts on the annulus $B_{2r}\setminus B_r$).
\end{remark}
The central estimate for the proof of our main result (Theorem \ref{MainResult}) reads as follows.
\begin{proposition}[Hole-filling estimate for 2D thin-film equation \eqref{tfe}]\label{LemmaHoldFill_TFE}
Let the assumptions of Theorem~\ref{MainResult} be satisfied. Let $b_r$ and $c_r$ be defined as per \eqref{defb}--\eqref{defc}. Then there exist constants $C_{\text{good}}>0, C_{\text{bad}}>0$, and $\delta\in (0,1/2]$ (all of these constants depending on $n$ only) such that, for any $r>0$, the estimate  
\begin{align}\label{hole_filling_total}
&\int_{B_{\delta r}} \frac12|\nabla u-b_r \cdot x -c_r|^2(\cdot,t_2) \d x
\\& \nonumber
\quad + C_{\text{\emph{good}}}\int_{t_1}^{t_2} \chi_{\text{\emph{good time}}(\delta r)}(t) \int_{B_{\delta r}} u^n |\nabla \Delta u|^2 + u^n |D^3 u|^2 \d x \d t
\\&\nonumber
\quad + C_{\text{\emph{bad}}}\int_{t_1}^{t_2} \chi_{\text{\emph{bad time}}(\delta r)}(t) \int_{B_{\delta r}} u^n |\nabla \Delta u|^2 + |\nabla u^{(n+2)/6}|^{6} \d x \d t
\\& \nonumber
\leq
\int_{B_{2r}} \frac12|\nabla u-b_r \cdot x -c_r|^2(\cdot,t_1) \d x
\\&~~~ \nonumber
+(1-\theta)\cdot C_{\text{\emph{good}}} \int_{t_1}^{t_2} \chi_{\text{\emph{good time}}(2r)}(t) \int_{B_{2 r}} u^n |\nabla \Delta u|^2 + u^n |D^3 u|^2 \d x \d t
\\&~~~ \nonumber
+ (1-\theta)\cdot C_{\text{\emph{bad}}} \int_{t_1}^{t_2} \chi_{\text{\emph{bad time}}(2r)}(t) \int_{B_{2 r}} u^n |\nabla \Delta u|^2 + |\nabla u^{(n+2)/6}|^{6} \d x \d t
\end{align}
holds with some $\theta=\theta(n)\in (0,1)$.
\end{proposition}
We recall the succinct, notationally convenient form of \eqref{hole_filling_total} that was introduced in \eqref{hf_est_succinct_tfe}.
\begin{proof}

The starting point is \eqref{123b} from Lemma~\ref{ExcessTimeEvolution}: This lemma allows to control the evolution of the tilt-excess-type quantity \eqref{abstract_def_tilt} via two specific terms. Both of these terms require different bounds for good and bad times. Therefore, we need four results, and these are given by Lemma~\ref{LemmaA}, Lemma~\ref{LemmaBC} (for good times), and Lemma~\ref{LemmaAdeg}, 
Lemma~\ref{LemmaBCdeg} (for bad times). Using these four lemmas in combination with 
\eqref{123b} gives, for a cutoff $\eta$  supported in $B_{2r}$ and equal to one in $B_{r}$ with $r^3|D^3 \eta|+r^2|D^2 \eta|+r|\nabla \eta|+\eta\leq C$,
\begin{align}\label{first_formula}
&\int_{B_{2r}} \frac12|\nabla u-b_r \cdot x -c_r|^2 \eta^6 \d x \bigg|_{t_1}^{t_2} + \int_{t_1}^{t_2} \int_{B_{2r}} u^n |\nabla \Delta u|^2\eta^6 \d x \d t
\\& \nonumber
\leq
C \int_{t_1}^{t_2} \chi_{\text{good time}(2r)}(t) \int_{B_{2r}\setminus B_{\delta r}} u^n |D^3 u|^2 \d x \d t
\\&~~~ \nonumber
+C \delta \int_{t_1}^{t_2} \chi_{\text{good time}(2r)}(t) \int_{B_{\delta r}} u^n |D^3 u|^2 \d x \d t
\\&~~~ \nonumber
+C \int_{t_1}^{t_2} \chi_{\text{bad time}(2r)}(t) \int_{B_{2r}\setminus B_{\delta r}} u^n |\nabla \Delta u|^2 + |\nabla u^{(n+2)/6}|^6 \d x \d t
\\&~~~ \nonumber
+C \delta \int_{t_1}^{t_2} \chi_{\text{bad time}(2r)}(t) \int_{B_{\delta r}} |\nabla u^{(n+2)/6}|^6 \d x \d t.
\end{align}
This estimate is still missing the $\int_{B_{\delta r}} u^n |D^3 u|^2 \d x$ term for good times and the $\int_{B_{\delta r}} |\nabla u^{(n+2)/6}|^6 \d x$ term for bad times on the left-hand side.

If $t$ is a good time for radius $2r$ (and thus, for any smaller radius as well, se Remark \ref{remark_good_times}), we can use Lemma~\ref{ThirdDerivativeBound} and write
\begin{align}\label{second_formula}
&\chi_{\text{good time}(2r)}(t) \int_{B_{r/2}} u^n |D^3 u|^2 \d x
\\& \nonumber
\leq C\chi_{\text{good time}(2r)}(t) \left\{\sup_{x\in B_{r}}{u^n}\!\cdot\! \int_{B_{r}} |\nabla \Delta u|^2 \d x
+C\!\! \sup_{x\in B_{r}}{u^n}\cdot\int_{B_{r}\setminus B_{r/2}} \!\!\!|D^3 u|^2 \d x\right\}
\\& \nonumber
\leq C\chi_{\text{good time}(2r)}(t) \left\{\inf_{x\in B_{r}}{u^n}\!\cdot\! \int_{B_{r}} |\nabla \Delta u|^2 \d x
+C\!\! \inf_{x\in B_r\setminus B_{r/2}}{\!\!\!\!\!u^n}\cdot\int_{B_{r}\setminus B_{r/2}} \!\!\!|D^3 u|^2 \d x\right\}
\\& \nonumber
\leq C\chi_{\text{good time}(2r)}(t) \left\{\int_{B_{r}} u^n |\nabla \Delta u|^2 \d x
+C \int_{B_{r}\setminus B_{r/2}} u^n |D^3 u|^2 \d x\right\}.
\end{align}

If $t$ is a bad time for $\delta r$, then $t$ is also a bad time for $2r$ (see Remark \ref{remark_good_times}). In this case, we use the Bernis-Gr\"un inequality \eqref{BernisGruenWeighted} with a cutoff $\eta$ with $\eta\equiv 1$ in $B_{r/2}$ and $\eta\equiv 0$ outside of $B_{r}$ to obtain
\begin{align}\label{third_formula}
&\chi_{\text{bad time}(\delta r)}(t) \int_{B_{r/2}} |\nabla u^{(n+2)/6}|^6 \d x
\\& \nonumber
\quad \leq C \chi_{\text{bad time}(\delta r)}(t) \int_{B_{r}} u^n |\nabla \Delta u|^2 + r^{-6} u^{n+2} \d x
\\& \quad \nonumber
\leq C \chi_{\text{bad time}(2 r)}(t) \int_{B_{r}} u^n |\nabla \Delta u|^2 \d x
\\&\quad \quad \nonumber
+C \chi_{\text{bad time}(2 r)}(t) \int_{B_{2r}\setminus B_{\delta r}} |\nabla u^{(n+2)/6}|^6 \d x
\\&\quad \quad \nonumber
+C \delta \chi_{\text{bad time}(2 r)}(t) \int_{B_{\delta r}} |\nabla u^{(n+2)/6}|^6 \d x, 
\end{align}
where in the last step we have used Lemma~\ref{BoundMaxUDegenerate} to estimate the integral $\int_{B_r}{u^{n+2}\d x}$, the fact that $\chi_{\text{bad time}(\delta r)}(t)\leq \chi_{\text{bad time}(2 r)}(t)$ (see Remark \ref{remark_good_times}),
and we have expanded the radius (from $r$ to $2r$) to get the second and third terms in the right-hand-side.

If $t$ is a bad time for $2r$ but a good time for $\delta r$, we use Lemma~\ref{LemmaThirdDerivativeByDissipation} to deduce
\begin{align}\label{fourth_formula}
&\chi_{\text{bad time}(2r)}(t) \chi_{\text{good time}(\delta r)}(t) \int_{B_r} u^n |D^3 u|^2 \d x
\\& \quad \nonumber
\leq
C \chi_{\text{bad time}(2r)}(t) \left\{\int_{B_{2r}} u^n |\nabla \Delta u|^2 \d x
+ \int_{B_{2r}} |\nabla u^{(n+2)/6}|^6 \d x\right\}. \nonumber
\end{align}
We now proceed to suitably combine \eqref{first_formula}, \eqref{second_formula}, \eqref{third_formula}, and \eqref{fourth_formula}. First, 
let now $\kappa \in(0,1]$ (we will impose conditions on it shortly). Summing up \eqref{second_formula} and \eqref{fourth_formula} weighted by $\kappa$, i.e., doing $\kappa \times \eqref{second_formula} + \kappa \times \eqref{fourth_formula}$, and taking into account that $\chi_{\text{good time}(2r)}(t) = \chi_{\text{good time}(2r)}(t)\chi_{\text{good time}(\delta r)}(t)$ (Remark \ref{remark_good_times}) and that $B_r\setminus B_{r/2}\subset B_{2r}\setminus B_{\delta r}$ (since $\delta \leq 1/2$) gives
\begin{align}\label{sum_kappa_eqns}
    & \kappa \left(\chi_{\text{good time}(\delta r)}(t) \int_{B_{r/2}} u^n |D^3 u|^2 \d x\right) \\
    &\nonumber \quad \leq C\kappa  \left(\int_{B_{2r}\setminus B_{\delta r}}{u^n |\nabla \Delta u|^2 \d x} + \int_{B_{\delta r}}{u^n |\nabla \Delta u|^2 \d x}\right)  \\
    &\nonumber \quad \quad + C\kappa \chi_{\text{good time}(2 r)}(t)\int_{B_{2r}\setminus B_{\delta r}}{u^n |D^3 u|^2 \d x} \\
    & \nonumber \quad \quad + C\kappa\chi_{\text{bad time}(2r)}(t) \int_{B_{2r}} |\nabla u^{(n+2)/6}|^6 \d x.
\end{align}

Performing the weighted expression sum $\eqref{first_formula} + \eqref{sum_kappa_eqns} + \tfrac{1}{2C}\times\eqref{third_formula}$
and using the fact that $\eta\equiv 1$ in $B_r$, we arrive at
\begin{align}\label{unabsorbed_tilt_excess_estimate}
 \sum_{i=1}^{4}{L_i} & := \int_{B_{2r}} \frac12|\nabla u-b_r \cdot x -c_r|^2 \eta^6 \d x \bigg|_{t_1}^{t_2}
\\&~~~ \nonumber
+ \int_{t_1}^{t_2} \int_{B_{r}} u^n |\nabla \Delta u|^2 \d x \d t
\\&~~~ \nonumber
+\kappa \int_{t_1}^{t_2} \chi_{\text{good time}(\delta r)}(t) \int_{B_{r/2}} u^n |D^3 u|^2 \d x \d t
\\&~~~ \nonumber
+\frac{1}{2C} \int_{t_1}^{t_2} \chi_{\text{bad time}(\delta r)}(t) \int_{B_{\delta r}} |\nabla u^{(n+2)/6}|^6 \d x \d t
\\& \nonumber
\leq
C \int_{t_1}^{t_2} \chi_{\text{good time}(2r)}(t) \int_{B_{2r}\setminus B_{\delta r}} u^n |\nabla \Delta u|^2 + u^n |D^3 u|^2 \d x \d t
\\&~~~ \nonumber
+C \delta \int_{t_1}^{t_2} \chi_{\text{good time}(2r)}(t) \int_{B_{\delta r}} u^n |D^3 u|^2 \d x \d t
\\&~~~ \nonumber
+C \int_{t_1}^{t_2} \chi_{\text{bad time}(2r)}(t) \int_{B_{2r}\setminus B_{\delta r}} u^n |\nabla \Delta u|^2 + |\nabla u^{(n+2)/6}|^6 \d x \d t
\\&~~~ \nonumber
+(C\delta + C \kappa) \int_{t_1}^{t_2} \chi_{\text{bad time}(2r)}(t) \int_{B_{\delta r}} |\nabla u^{(n+2)/6}|^6 \d x \d t
\\&~~~ \nonumber
+C \kappa \int_{t_1}^{t_2} \int_{B_{\delta r}} u^n |\nabla \Delta u|^2 \d x \d t 
\\&~~~ \nonumber 
+ \frac{1}{2}\int_{t_1}^{t_2}\chi_{\text{bad time}(2 r)}(t) \int_{B_{r}} u^n |\nabla \Delta u|^2 \d x \d t =: \sum_{i=1}^{6}{R_i},
\end{align}
where we also renamed $C$ (on the account of having $\kappa\leq 1$). 

We now work on the right-hand-side of \eqref{unabsorbed_tilt_excess_estimate} to be able to perform hole-filling. We keep the terms $R_1$ are $R_3$ as they are, as these are ready for hole-filling.
For $\kappa$ small enough, we can absorb $R_5+R_6$ into $L_2$. Furthermore, for $0<\delta\ll \kappa$ small enough, and noticing that $\chi_{\text{good time}(2 r)}(t) \leq \chi_{\text{good time}(\delta r)}(t)$, we can absorb $R_2$. The only term that is left is $R_4$, which can not be absorbed (because, if $t$ is a bad time for $2r$, $t$ may not necessarily be bad for $\delta r$). However, the term $R_4$ has the same form as the term that will be added during the hole-filling procedure and -- crucially -- is multiplied with the arbitrarily small factor $(\delta + \kappa)$, so it will not create problems once $R_1$ and $R_3$ get hole-filled.
Performing the absorptions in \eqref{unabsorbed_tilt_excess_estimate} as previously detailed, we obtain, for some constant $\tilde C=\tilde C(n)>0$,
\begin{align*}
&\int_{B_{2r}} \frac12|\nabla u-b_r \cdot x -c_r|^2 \eta^6 \d x \bigg|_{t_1}^{t_2}
\\& \quad 
+\kappa \int_{t_1}^{t_2} \chi_{\text{good time}(\delta r)}(t) \int_{B_{\delta r}} u^n |\nabla \Delta u|^2 + u^n |D^3 u|^2 \d x \d t
\\& \quad 
+\frac{1}{\tilde C} \int_{t_1}^{t_2} \chi_{\text{bad time}(\delta r)}(t) \int_{B_{\delta r}} u^n |\nabla \Delta u|^2 + |\nabla u^{(n+2)/6}|^6 \d x \d t
\\&
\leq
\tilde C \int_{t_1}^{t_2} \chi_{\text{good time}(2r)}(t) \int_{B_{2r}\setminus B_{\delta r}} u^n |\nabla \Delta u|^2 + u^n |D^3 u|^2 \d x \d t
\\&~~~
+\tilde C \int_{t_1}^{t_2} \chi_{\text{bad time}(2r)}(t) \int_{B_{2r}\setminus B_{\delta r}} u^n |\nabla \Delta u|^2 + |\nabla u^{(n+2)/6}|^6 \d x \d t
\\&~~~
+ (\tilde{C}\delta + \tilde{C}\kappa) \int_{t_1}^{t_2} \chi_{\text{bad time}(2r)}(t) \int_{B_{\delta r}} |\nabla u^{(n+2)/6}|^6 \d x \d t.
\end{align*}
This entails, by hole-filling, that
\begin{align}\label{once_hole_filled}
&\int_{B_{2r}} \frac12|\nabla u-b_r \cdot x -c_r|^2 \eta^6 \d x \bigg|_{t_1}^{t_2}
\\&\nonumber \quad
+(\tilde C+\kappa) \int_{t_1}^{t_2} \chi_{\text{good time}(\delta r)}(t) \int_{B_{\delta r}} u^n |\nabla \Delta u|^2 + u^n |D^3 u|^2 \d x \d t
\\&\nonumber \quad
+(\tilde C+\tfrac{1}{\tilde C}) \int_{t_1}^{t_2} \chi_{\text{bad time}(\delta r)}(t) \int_{B_{\delta r}} u^n |\nabla \Delta u|^2 + |\nabla u^{(n+2)/6}|^6 \d x \d t
\\&\nonumber
\leq
\tilde C \int_{t_1}^{t_2} \chi_{\text{good time}(2r)}(t) \int_{B_{2r}} u^n |\nabla \Delta u|^2 + u^n |D^3 u|^2 \d x \d t
\\&~~~\nonumber
+\tilde{C}(1+\kappa+\delta) \int_{t_1}^{t_2} \chi_{\text{bad time}(2r)}(t) \int_{B_{2r}} u^n |\nabla \Delta u|^2 + |\nabla u^{(n+2)/6}|^6 \d x \d t.
\end{align}
For $\kappa$ sufficiently small (which also implies $\delta$ sufficiently small by the previously imposed relation $\delta \ll \kappa$), it holds 
$$
1-\theta := \max\left\{\tfrac{\tilde{C}}{\tilde{C}+\kappa};\tfrac{\tilde{C}(1+\kappa+\delta)}{\tilde{C}+\tfrac{1}{\tilde{C}}}\right\} \in (0,1).
$$
Therefore, using the fact that $\eta \equiv 1$ on $B_{\delta r},$ we see that \eqref{once_hole_filled} entails our desired estimate \eqref{hole_filling_total} upon setting $C_{\text{good}} := \tilde{C}+\kappa$ and $C_{\text{bad}} := \tilde{C}+\tfrac{1}{\tilde{C}}$.
\end{proof}

\subsection{An evolution equation for the tilt-excess-type quantity \eqref{abstract_def_tilt}}

In the next lemma, we provide the basic estimate on the evolution of the tilt-excess-type quantity \eqref{abstract_def_tilt} (via a close analogue which enjoys smoothing via a compactly supported test function $\eta$) for the solution to the thin-film equation.

\begin{lemma}
\label{ExcessTimeEvolution}
Let the assumptions of Theorem~\ref{MainResult} be satisfied.
Let $r>0$ and let $\eta$ be a cutoff supported in $B_{2r}$ and equal to one in $B_{r}$ with $r^3|D^3 \eta|+r^2|D^2 \eta|+r|\nabla \eta|+\eta\leq C$.
Then the solution $u$ admits the estimate
\begin{align}\label{123b}
&\int_{B_{2r}} \frac12|\nabla u-b_r \cdot x -c_r|^2 \eta^6 \d x \bigg|_{t_1}^{t_2} + \int_{t_1}^{t_2} \int_{B_{2r}} u^n |\nabla \Delta u|^2\eta^6 \d x \d t
\\&
\nonumber
\leq
-\int_{t_1}^{t_2} \int_{B_{2r}} u^n \nabla \Delta u \cdot \Big\{ (\Delta u-\operatorname{tr} b_r) \nabla \eta^6  \\
& \nonumber \quad \quad \quad \quad \quad \quad + (D^2u-b_r) \cdot \nabla \eta^6 + (\nabla u -b_r \cdot x -c_r)\cdot D^2 \eta^6 \Big\} \d x \d t
\\&~~~
\nonumber
- \int_{t_1}^{t_2} \int_{B_{2r}}(b_r'\cdot x + c_r')\cdot (\nabla u -b_r \cdot x -c_r) \eta^6 \d x \d t.
\end{align}
\end{lemma}
To establish \eqref{123b}, we make use of the following weighted energy dissipation principle proved in \cite{DeNittiFischer} for energy-dissipating weak solutions to the thin-film equation.
\begin{lemma}[Weighted energy estimate, see Lemma A.3 in \cite{DeNittiFischer}]
	\label{WeightedEnergy}
	Let $\Omega = \R^d$, let $n \in (2-\sqrt{4/5}, 3)$,  and let $u$ be an  energy-dissipating weak solution to the thin-film equation \eqref{tfe} with zero contact angle in the sense of Definition~\ref{DefinitionEnergyDissipatingWeakSolution}. Let $\psi\in
	C^2_{cpt}(\R^d)$ be a nonnegative weight function. Then we have
	\begin{align}\label{EnergyEstimateA}
	&\int_{\R^d} \frac{1}{2}|\nabla u|^2\psi \d x\Bigg|_{t_1}^{t_2}
	-\int_{t_1}^{t_2}\int_{\R^d} \frac{1}{2}|\nabla u|^2 \partial_t \psi \d x\d t
	\\ \nonumber
    &~~~=-\int_{t_1}^{t_2}\int_{\{u(\cdot, t)>0\}} u^n|\nabla \Delta u|^2\psi \d x\d t
	\\&
	\qquad \nonumber-\int_{t_1}^{t_2}\int_{\{u(\cdot, t)>0\}}u^n\nabla\Delta u\cdot
	\Big(\Delta u\nabla \psi+D^2u\cdot\nabla \psi+\nabla u\cdot D^2\psi\Big)
	\d x \d t
	\end{align}
	for a.e. $t_2\geq t_1\geq 0$ and a.e. $t_2\geq 0$ in case $t_1=0$.
\end{lemma}

\begin{proof}[Proof of Lemma~\ref{ExcessTimeEvolution}]
We expand 
\begin{align}\label{123}
\int_{B_{2r}} &\frac12|\nabla u-b_r \cdot x -c_r|^2 \eta^6 \d x
\bigg|_{t_1}^{t_2}
\\&\nonumber
=\int_{B_{2r}} \frac12 |\nabla u|^2\eta^6 \d x 
\bigg|_{t_1}^{t_2}
-\int_{B_{2r}} \nabla u \cdot (b_r \cdot x +c_r)\eta^6 \d x
\bigg|_{t_1}^{t_2} \\
& \quad \nonumber
+ \int_{B_{2r}} \frac12|b_r \cdot x +c_r|^2 \eta^6 \d x
\bigg|_{t_1}^{t_2}.
\end{align}
We treat the three integrals on the right-hand-side of \eqref{123} separately. For the first one, we apply Lemma~\ref{WeightedEnergy}: In \eqref{EnergyEstimateA}, we substitute $\psi:=\eta^6$ and use that it is time-independent and that $r^3|D^3 \psi|+r^2|D^2 \psi|+r|\nabla \psi|+\eta\leq C$ to obtain
\begin{align}\label{1}
&\int_{B_{2r}} \frac12 |\nabla u|^2\eta^6 \d x
\bigg|_{t_1}^{t_2}
\\ \nonumber
&= - \int_{t_1}^{t_2} \int_{B_{2r}} u^n |\nabla \Delta u|^2\eta^6 \d x \d t \\
&~~~\nonumber
-\int_{t_1}^{t_2} \int_{B_{2r}} u^n \nabla \Delta u \cdot \left( \Delta u \nabla \eta^6 + D^2 u \cdot \nabla \eta^6 + \nabla u \cdot D^2 \eta^6 \right) \d x \d t.
\end{align}
For the second integral on the right-hand side of \eqref{123}, we apply~\eqref{def:endisWS} from Definition~\ref{DefinitionEnergyDissipatingWeakSolution} with $\psi:=\nabla \cdot ((b_r \cdot x +c_r)\eta^6)$ and rely on the symmetry of $b_r$,
thus obtaining (also abbreviating $\langle \cdot, \cdot \rangle = \langle \cdot, \cdot\rangle_{(W^{1,p}(\R^d))'\times W^{1,p}(\R^d)}$ for some $p > \frac{4d}{2d + n(2-d)}$)
\begin{align}\label{2}
&-\int_{B_{2r}} \nabla u \cdot (b_r \cdot x +c_r)\eta^6 \d x \bigg|_{t_1}^{t_2}
\\& \nonumber
\quad = \int_{t_1}^{t_2} \langle \partial_t u, \nabla \cdot ((b_r \cdot x +c_r)\eta^6) \rangle \d t
- \int_{t_1}^{t_2} \int_{B_{2r}} \nabla u \cdot (b_r' \cdot x + c_r') \eta^6 \d x \d t
\\
& \nonumber \quad \stackrel{\eqref{def:endisWS}}{=} \int_{t_1}^{t_2}\int_{B_{2r}} u^n \nabla \Delta u \cdot \nabla \left( \nabla \cdot ((b_r \cdot x +c_r)\eta^6) \right) \d x \d t
\\& \quad ~~~ \nonumber
-\int_{t_1}^{t_2} \int_{B_{2r}} \nabla u \cdot (b_r' \cdot x + c_r') \eta^6 \d x \d t
\\& \nonumber
 \quad  =\int_{t_1}^{t_2} \int_{B_{2r}} u^n \nabla \Delta u \cdot \left(\operatorname{tr} b \, \nabla \eta^6 + b_r \cdot \nabla \eta^6 + (b_r \cdot x +c_r)\cdot D^2\eta^6 \right) \d x \d t
\\& \quad ~~~ \nonumber
-
\int_{t_1}^{t_2} \int_{B_{2r}} \nabla u \cdot (b_r' \cdot x + c_r') \eta^6 \d x \d t.
\end{align}
Finally, we simply rewrite the last integral in \eqref{123} as 
\begin{equation}\label{3}
\int_{B_{2r}} \frac12 |b_r \cdot x +c_r|^2 \eta^6 \d x \bigg|_{t_1}^{t_2} = \int_{t_1}^{t_2} \int_{B_{2r}}(b_r'\cdot x + c_r')\cdot (b_r \cdot x +c_r) \eta^6 \d x \d t.
\end{equation}
Altogether, using~\eqref{1},~\eqref{2} and~\eqref{3} in~\eqref{123} we obtain \eqref{123b}.
\end{proof}
We now need to estimate the two terms appearing in the right-hand-side of \eqref{123b}. For doing this, we treat separately the case of good times (Subsection \ref{subsec_good}) and bad times (Subsection \ref{subsec_bad}). Furthermore, when transitioning from a bad time to a good time (due to the decrease in $r$), we need an argument to control $\int u^n |D^3 u|^2 \d x$ on the smaller ball; this is provided in Subsection \ref{subsec_transition}.

\subsection{The estimate for good times}\label{subsec_good}


For a good time $t$ of the weak solution to the thin-film equation (see Definitions \ref{DefinitionEnergyDissipatingWeakSolution} and \ref{DefnGoodBad}), we estimate the terms on the right hand side of~\eqref{123b} separately in the following two lemmas.


\begin{lemma}\label{LemmaA}
Let the assumptions of Theorem~\ref{MainResult} be satisfied. Let $r>0$ and let $\eta$ be as in Lemma~\ref{ExcessTimeEvolution}.
If $t$ is a good time for the radius $2r$ (in the sense of Definition~\ref{DefnGoodBad}), we can estimate
\begin{equation}\label{eqlemmaA}
\begin{aligned}
-\int_{B_{2r}} \!\!\!\! &u^n \nabla \Delta u \cdot \left( (\Delta u-\operatorname{tr} b_r) \nabla \eta^6 + (D^2u-b_r) \cdot \nabla \eta^6 + (\nabla u -b_r \cdot x -c_r)\cdot D^2 \eta^6 \right) \!\! \d x
\\
&\leq C \Big(\inf_{B_{2r}} u\Big)^n \int_{B_{2r}\setminus B_r} |D^3 u|^2 \d x.
\end{aligned}
\end{equation}
\end{lemma}
\begin{proof}
We use the properties of $\eta$ and Young's inequality to estimate
\begin{align*}
&\mbox{LHS of }\eqref{eqlemmaA}\\ 
& \quad \leq C \Big(\sup_{B_{2r}} u\Big)^n \int_{B_{2r}\setminus B_r} |D^3 u|\left(|D^2 u-b_r| r^{-1} + |\nabla u-b_r \cdot x -c_r| r^{-2}\right) \d x
\\
& \quad \leq C \Big(\sup_{B_{2r}} u\Big)^n \int_{B_{2r}\setminus B_r} |D^3 u|^2 \d x + C  \Big(\sup_{B_{2r}} u\Big)^n r^{-2}\int_{B_{2r}\setminus B_r} |D^2 u-b_r|^2 \d x
\\& \quad \quad 
+ C \Big(\sup_{B_{2r}} u\Big)^n r^{-4}\int_{B_{2r}\setminus B_r}|\nabla u-b_r \cdot x -c_r|^2 \d x.
\end{align*}
Using the Poincar\'e inequality on the annulus $B_{2r}\setminus B_r$ (see Lemma \ref{lemma_poincare}), we infer
\begin{align*}
\int_{B_{2r}\setminus B_r}|\nabla u-b_r \cdot x -c_r|^2 \d x
& \leq C r^2 \int_{B_{2r}\setminus B_r}|D^2 u-b_r|^2 \d x
\\&
\leq C r^4 \int_{B_{2r}\setminus B_r}|D^3 u|^2 \d x.
\end{align*}
Combining the previous two estimates with the definition of a good time \eqref{goodt} concludes the proof.
\end{proof}

\begin{lemma}\label{LemmaBC}
Let the assumptions of Theorem~\ref{MainResult} be satisfied.  Let $r>0$ and let $\eta$ be as in Lemma~\ref{ExcessTimeEvolution}.
If $t$ is a good time for the radius $2r$ (in the sense of Definition~\ref{DefnGoodBad}), we can estimate for any $\delta\in (0,1]$
\begin{align}
\label{eqlemmaBC}
&\left| \int_{B_{2r}}(b_r'\cdot x + c_r')(\nabla u-b_r \cdot x -c_r) \eta^6 \d x \right| 
\\&\nonumber
\quad \leq C \Big(\inf_{B_{2r}} u\Big)^n \int_{B_{2r}\setminus B_{\delta r}} |D^3 u|^2 \d x
+C \delta \Big(\inf_{B_{2r}} u\Big)^n \int_{B_{\delta r}} |D^3 u|^2 \d x.
\end{align}
\end{lemma}
\begin{proof}
From~\eqref{defb} and~\eqref{defc}, using Definition~\ref{DefinitionEnergyDissipatingWeakSolution}, in particular~\eqref{def:endisWS}, we can compute
\begin{align}
b_r'(t)&= \frac{1}{\int_{B_{2r}} \tilde{\eta}(\frac{x}{r})\d x} \int_{B_{2r}} D^3 \Big( \tilde{\eta}\left(\frac{x}{r}\right) \Big) \cdot u^n \nabla \Delta u \d x, \label{defb'}
\\
c_r'(t)&=-\frac{1}{\int_{B_{2r}} \tilde{\eta}(\frac{x}{r})\d x} \int_{B_{2r}} D^2 \Big( \tilde{\eta}\left(\frac{x}{r}\right) \Big) \cdot u^n \nabla \Delta u \d x.\label{defc'}
\end{align}
This entails by H\"older's inequality, the bound $r^3|D^3 \tilde \eta| + r^2|D^2 \tilde \eta| + r|\nabla \tilde \eta| + \tilde \eta \leq C$, the fact that $\supp \tilde \eta \subset B_{2r} \setminus B_r$, as well as the definition of a good time \eqref{goodt}, that
\begin{align}
\label{BoundB'C'}
r |b_r'|+|c_r'|
\leq C r^{-d/2-2} \Big(\inf_{B_{2r}} u\Big)^n \bigg(\int_{B_{2r}\setminus B_{r}} |D^3 u|^2 \d x\bigg)^{1/2}.
\end{align}
Furthermore, using the Poincar\'e inequality and the Poincar\'e-Sobolev inequality (both on $B_{2r}$, see Lemma \ref{lemma_poincare}) entail
\begin{align*}
&\int_{B_{2r}} |\nabla u-b_r \cdot x -c_r|^2 \d x
\\& \quad 
\leq C r^2 \int_{B_{2r}} |D^2 u-b_r|^2 \d x
\\& \quad
\leq C r^{4-d} \bigg(\int_{B_{2r}} |D^3 u| \d x\bigg)^2
\\& \quad
\leq C r^{4-d} \bigg(\int_{B_{2r}\setminus B_{\delta r}} |D^3 u| \d x\bigg)^2 +C r^{4-d} \bigg(\int_{B_{\delta r}} |D^3 u| \d x \bigg)^2.
\end{align*}
Applying H\"older's inequality, we deduce
\begin{align*}
&\int_{B_{2r}} |\nabla u-b_r \cdot x -c_r|^2 \d x
\leq
C r^{4} \int_{B_{2r}\setminus B_{\delta r}} \!\! |D^3 u|^2 \d x
+ C \delta^d r^{4} \int_{B_{\delta r}} |D^3 u|^2 \d x.
\end{align*}
Combining this estimate with \eqref{BoundB'C'}, our lemma follows.
\end{proof}

\begin{lemma}
\label{ThirdDerivativeBound}
Let $u\in H^3(B_{2r})$.
We then have the estimate
\begin{align}\label{ineq_in_ThirdDerivativeBound}
\int_{B_r} |D^3 u|^2 \d x
\leq C\int_{B_{2r}} |\nabla \Delta u|^2 \d x
+C\int_{B_{2r}\setminus B_r} |D^3 u|^2 \d x.
\end{align}
\end{lemma}
\begin{proof}
Let $u$ be a smooth function and let $r>0$. Set $\bar b:=\dashint_{B_{2r}\setminus B_r} D^2 u \d x$. Let $\eta$ be a cutoff supported in $B_{2r}$ with $\eta\equiv 1$ in $B_r$. We then have
\begin{align*}
&\int |D^3 u|^2 \eta^2 \d x
\\
& \quad
=-\int (D^2 u-\bar b) : D^2 \Delta u ~ \eta^2 \d x
-2\int \eta (D^2 u-\bar b) : D^3 u \cdot \nabla \eta \d x
\\& \quad 
=\int |\nabla \Delta u|^2 \eta^2 \d x
+2\int \eta \nabla \Delta u \cdot (D^2 u-\bar b) \cdot \nabla \eta \d x \\
& \quad \quad -2\int \eta (D^2 u-\bar b) : D^3 u \cdot \nabla \eta \d x.
\end{align*}
where we also added $-\overline{b}$ in the first equality.
Using Young's inequality and absorption, we arrive at
\begin{align*}
\int |D^3 u|^2 \eta^2 \d x
& \leq C \int |\nabla \Delta u|^2 \eta^2 \d x
+ C \int |D^2 u-\bar b|^2 |\nabla \eta|^2 \d x
\\&
\leq C \int |\nabla \Delta u|^2 \eta^2 \d x
+ C r^{-2} \int_{B_{2r}\setminus B_r} |D^2 u-\bar b|^2 \d x.
\end{align*}
By approximation, this estimate holds for any $u\in H^3(B_{2r})$.
The Poincar\'e inequality on the annulus $B_{2r}\setminus B_r$ now yields the desired estimate.
\end{proof}

\subsection{Transitioning from a bad time to a good time}\label{subsec_transition}

\begin{lemma}
\label{LemmaThirdDerivativeByDissipation}
Let the assumptions of Theorem~\ref{MainResult} be satisfied.  Let $r>0$ and let $\eta$ be as in Lemma~\ref{ExcessTimeEvolution}.
If $t$ is a bad time for $2r$ but a good time for $\delta r$, we have 
\begin{align*}
\int_{B_r} u^n |D^3 u|^2 \d x
\leq C \int_{B_{2r}} u^n |\nabla \Delta u|^2 \d x
+C \int_{B_{2r}} |\nabla u^{(n+2)/6}|^6 \d x.
\end{align*}
\end{lemma}
\begin{proof}
Since the cutoff $\eta$ is as in Lemma \ref{ExcessTimeEvolution}, it satisfies $\eta\equiv 1$ in $B_r$ and $\eta \equiv 0$ outside of $B_{2r}$. Using the elementary formula
\begin{align*}
D^3 u^{(n+2)/2}&=\tfrac{n+2}{2} u^{n/2} D^3 u+3\tfrac{(n+2)n}{4} u^{(n-2)/2} D^2 u \otimes \nabla u
\\&~~~
+ \tfrac{(n+2)n(n-2)}{8} u^{(n-4)/2} \nabla u \otimes \nabla u \otimes \nabla u,
\end{align*}
the Bernis-Gr\"un formula \eqref{BernisGruenWeighted}, and the fact that $\nabla u^{(n+2)/6} = \tfrac{n+2}{6} u^{(n-4)/6}\nabla u$, we obtain
\begin{align}\label{control_un_full_third_derivative}
&\int_{B_r} u^n |D^3 u|^2 \d x
\leq C(d,n) \int_{B_{2r}} \!\! u^n |\nabla \Delta u|^2 \eta^6 \d x + C(d,n) r^{-6} \int_{B_{2r}\setminus B_r}\!\!\!\!\! u^{n+2} \d x.
\end{align}
Applying Lemma~\ref{BoundMaxUDegenerate} (with any $\delta\leq \tfrac{1}{2}$) to estimate the last term yields the result.
\end{proof}

\subsection{The estimate for bad times}\label{subsec_bad}

\begin{lemma}
\label{BoundMaxUDegenerate}
Let $n>0$ and $d=2$.
Let $u\in L^1(\mathbb{R}^2)$ be nonnegative and satisfy $\nabla u^{(n+2)/6}\in L^6(\mathbb{R}^2)$. Assume that $\min_{B_{2r}} u < \frac{1}{2} \max_{B_{2r}} u$.
For any $\delta \in (0,1]$ we then have the bound
\begin{align*}
\max_{B_{2r}} u^{n+2}
\leq C r^{6-d} \int_{B_{2r}\setminus B_{\delta r}} |\nabla u^{(n+2)/6}|^6 \d x + Cr^{6-d} \delta \int_{B_{\delta r}} |\nabla u^{(n+2)/6}|^6 \d x.
\end{align*}
\end{lemma}
\begin{proof}
By Morrey's inequality
(applicable since $3>d=2$, see Lemma~\ref{lemma_morrey_ineq}) and invariance of the quantity $\max_{B_{2r}} u^{(n+2)/6} - \min_{B_{2r}} u^{(n+2)/6}$ with respect to adding constants to the function $u^{(n+2)/6}$, we have
\begin{align*}
\max_{B_{2r}} u^{(n+2)/6} - \min_{B_{2r}} u^{(n+2)/6}
\leq
C r \bigg(\dashint_{B_{2r}} |\nabla u^{(n+2)/6}|^3 \d x \bigg)^{1/3}.
\end{align*}
Using the fact that $\min_{B_{2r}} u\leq \frac{1}{2} \max_{B_{2r}} u$ and raising both sides to the power $6$, we obtain 
\begin{align*}
\max_{B_{2r}} u^{n+2}
\leq
C(n) r^6 \bigg(r^{-d} \int_{B_{2r}} |\nabla u^{(n+2)/6}|^3 \d x \bigg)^{2}.
\end{align*}
Splitting the integral yields
\begin{align*}
\max_{B_{2r}} u^{n+2}
& \leq 
C(n) r^{6-2d} \bigg(\int_{B_{2r}\setminus B_{\delta r}} |\nabla u^{(n+2)/6}|^3 \d x \bigg)^2
\\&
 \quad + C(n) r^{6-2d} \bigg(\int_{B_{\delta r}} |\nabla u^{(n+2)/6}|^3 \d x \bigg)^{2}.
\end{align*}
An application of the H\"older inequality now yields the statement of the lemma.
\end{proof}

\begin{lemma}
\label{LemmaAdeg}
Let the assumptions of Theorem~\ref{MainResult} be satisfied. Let $r>0$ and let $\eta$ be as in Lemma~\ref{ExcessTimeEvolution}.
If $t$ is a bad time for the radius $2r$ (in the sense of Definition~\ref{DefnGoodBad}), we have the estimate
\begin{align}\label{TiltExcessBad}
&-\int_{B_{2r}} \!\!\! u^n \nabla \Delta u \cdot \left( (\Delta u-\operatorname{tr} b_r) \nabla \eta^6 + (D^2u-b_r) \cdot \nabla \eta^6 + (\nabla u -b_r \cdot x -c_r)\cdot D^2 \eta^6 \right) \!\!\d x
\\&
\quad \leq
C \int_{B_{2r}\setminus B_r} u^n |\nabla \Delta u|^2 \d x
+ C  \int_{B_{2r} \setminus B_{\delta r}} |\nabla u^{(n+2)/6}|^6 \d x \nonumber
\\& \quad \quad
+ C \delta \int_{B_{\delta r}} |\nabla u^{(n+2)/6}|^6 \d x.\nonumber
\end{align}
\end{lemma}
\begin{proof}
Using Young's inequality and the fact that $\supp \nabla \eta \cup \supp D^2 \eta \subset B_{2r}\setminus B_r$ as well as $|\eta|\leq 1$ and $r^2 |D^2 \eta| + r |\nabla \eta|\leq C$, we obtain
\begin{align*}
\textrm{LHS of }\eqref{TiltExcessBad} &
\leq 
\int_{B_{2r}\setminus B_r} u^n |\nabla \Delta u|^2 \,\d x
\\&~~~
+C \int_{B_{2r}\setminus B_r} u^n (|D^2 u|^2 + |b_r|^2) |\nabla \eta|^2 \,\d x
\\&~~~
+C r^{-4} \int_{B_{2r}\setminus B_r} u^n (|\nabla u|^2 + r^2 |b_r|^2 + |c_r|^2) \,\d x.
\end{align*}
Using Lemma~\ref{SecondDerivativeEstimate} with $\eta$ replaced by $|\nabla \eta|$ to estimate the second term on the right-hand-side above and inserting the bound $r |\nabla \eta|\leq C$, we obtain
\begin{align*}
\textrm{LHS of }\eqref{TiltExcessBad}&
\leq
2\int_{B_{2r}\setminus B_r} u^n |\nabla \Delta u|^2 \,\d x
\\&~~~
+C \int u^n |\nabla u|^2 (|\nabla \eta|^4+|D^2 \eta|^2) \d x
\\&~~~
+C \int u^{n-2} |\nabla u|^4 |\nabla \eta|^2 \d x
\\&~~~
+C r^{-4} \int_{B_{2r}\setminus B_r} u^n |\nabla u|^2 \,\d x
\\&~~~
+C r^{d-4} (r^2 |b_r|^2 + |c_r|^2) \max_{B_{2r}\setminus B_r} u^n .
\end{align*}
Using again $r^2 |D^2 \eta| + r |\nabla \eta|\leq C$, we deduce
\begin{align*}
\textrm{LHS of }\eqref{TiltExcessBad}&
\leq
2\int_{B_{2r}\setminus B_r} u^n |\nabla \Delta u|^2 \,\d x
\\&~~~
+C r^{-4} \max_{B_{2r}\setminus B_r} u^{4(n+2)/6} \int_{B_{2r}\setminus B_r} |\nabla u^{(n+2)/6}|^2  \d x
\\&~~~
+C r^{-2} \max_{B_{2r}\setminus B_r} u^{2(n+2)/6} \int_{B_{2r}\setminus B_r} |\nabla u^{(n+2)/6}|^4 \d x
\\&~~~
+C r^{d-4} (r^2 |b_r|^2 + |c_r|^2) \max_{B_{2r}\setminus B_r} u^n .
\end{align*}
Using the square of the estimate
\begin{align}\label{rb_plus_c_estimate}
r^2 |b_r| + r |c_r| \leq C \max_{B_{2r}\setminus B_r} u
\end{align}
(estimate \eqref{rb_plus_c_estimate} simply follows by inspecting the definition of $b_r$ and $c_r$ given in \eqref{defb}--\eqref{defc}) and using the H\"older inequality, we get
\begin{align*}
\textrm{LHS of }\eqref{TiltExcessBad}&
\leq
2\int_{B_{2r}\setminus B_r} u^n |\nabla \Delta u|^2 \,\d x
\\&~~~
+C r^{2d/3-4} \max_{B_{2r}\setminus B_r} u^{4(n+2)/6} \bigg(\int_{B_{2r}\setminus B_r} |\nabla u^{(n+2)/6}|^6  \d x \bigg)^{1/3}
\\&~~~
+C r^{d/3-2} \max_{B_{2r}\setminus B_r} u^{2(n+2)/6} \bigg(\int_{B_{2r}\setminus B_r} |\nabla u^{(n+2)/6}|^6 \d x\bigg)^{2/3}
\\&~~~
+C r^{d-6} \max_{B_{2r}\setminus B_r} u^{n+2} .
\end{align*}
Inserting the bound from Lemma~\ref{BoundMaxUDegenerate} and using Young's inequality, we deduce our desired estimate.
\end{proof}

\begin{lemma}
\label{SecondDerivativeEstimate}
Let $\eta$ be a nonnegative Lipschitz function with compact support.
For any nonnegative function $u\in L^1(\mathbb{R}^2)$ with $\nabla u^{(n+2)/6}\in L^6(\mathbb{R}^2)$ and $u^{n/2} \nabla \Delta u \in L^2(\mathbb{R}^2)$, we have $D^2 u\in L^2_{loc}(\{x: u(x) \eta(x)>0\})$ and
\begin{align*}
\int u^n |D^2 u|^2 \eta^2 \d x
&\leq \int_{\supp \eta} u^n |\nabla \Delta u|^2 \d x
\\&~~~
+C \int u^n |\nabla u|^2 (\eta^4+|\nabla \eta|^2) \d x
\\&~~~
+C \int u^{n-2} |\nabla u|^4 \eta^2 \d x.
\end{align*}
\end{lemma}
\begin{proof}
For a nonnegative smooth function $u$, we integrate by parts to obtain
\begin{align*}
\int u^n |D^2 u|^2 \eta^2 \d x
&=-\int u^n \nabla \Delta u \cdot \nabla u ~ \eta^2 \d x
\\&~~~
-n \int u^{n-1} \nabla u \cdot D^2 u \cdot \nabla u ~ \eta^2 \d x
\\&~~~
-2\int u^n \nabla u \cdot D^2 u \cdot \eta \nabla \eta \d x.
\end{align*}
Applying Young's inequality and absorbing then yields the claim for nonnegative smooth functions $u$. The approximation of $u$ by smooth functions $u_\varepsilon$ given in the proof of \cite[Lemma~A.4]{DeNittiFischer} then justifies the estimate under the stated regularity requirements. Note that by the previous estimate applied to the approximations $u_\varepsilon$ and by the continuity of $u$, the restrictions $D^2 u_\varepsilon|_{U_\kappa}$ to any of the open sets $U_\kappa:=\{x:u(x)>\kappa,\eta(x)>\kappa\}$ converge weakly to $D^2 u|_{U_\kappa}$ in $H^2(U_\kappa)$, thereby implying the regularity $D^2 u\in L^2_{loc}(\{x:u(x)\eta(x)>0\})$.
\end{proof}

\begin{lemma}
\label{LemmaBCdeg}
Let the assumptions of Theorem~\ref{MainResult} be satisfied.  Let $r>0$ and let $\eta$ be as in Lemma~\ref{ExcessTimeEvolution}.
If $t$ is a bad time for the radius $2r$ (in the sense of Definition~\ref{DefnGoodBad}), we have for any $\delta\in (0,\tfrac{1}{2}]$
\begin{align*}
& \left| \int_{B_{2r}}(b_r'\cdot x + c_r')(\nabla u-b_r \cdot x -c_r) \eta^6 \d x \right|
\\& \quad
\leq
C \int_{B_{2r}\setminus B_r} u^n |\nabla \Delta u|^2 \d x
+ C \int_{B_{2r} \setminus B_{\delta r}} |\nabla u^{(n+2)/6}|^6 \d x
\\& \quad \quad 
+C \delta \int_{B_{\delta r}} |\nabla u^{(n+2)/6}|^6 \d x.
\end{align*}
\end{lemma}
\begin{proof}
We have
\begin{align*}
& \left|\int_{B_{2r}}(b_r'\cdot x + c_r')(\nabla u-b_r \cdot x -c_r) \eta^6 \d x \right| 
\\& \quad 
\leq
C (r|b_r'| + |c_r'|) \bigg(\int_{B_{2r}} |\nabla u| \d x + r^{d+1} |b_r| + r^d |c_r| \bigg).\nonumber
\end{align*}
A simple inspection of \eqref{defb'}--\eqref{defc'} grants
\begin{align}\label{control_der_b_c}
r^2|b_r'|+r|c_r'|
& 
\leq C r^{-d-1} \int_{B_{2r}\setminus B_r} u^n |\nabla \Delta u| \d x
\\& 
\leq C r^{-d/2-1} \max_{B_{2r}\setminus B_r} u^{n/2} \bigg(\int_{B_{2r}\setminus B_r} u^n |\nabla \Delta u|^2 \d x\bigg)^{1/2}.\nonumber
\end{align}
Using \eqref{control_der_b_c} with \eqref{rb_plus_c_estimate} entails
\begin{align*}
&\left| \int_{B_{2r}}(b_r'\cdot x + c_r')(\nabla u-b_r \cdot x -c_r) \eta^6 \d x \right|
\\& \quad
\leq
Cr^{-d/2-2} \max_{B_{2r}\setminus B_r} u^{n/2} \bigg(\int_{B_{2r}\setminus B_r} u^{n} |\nabla \Delta u|^2 \d x\bigg)^{1/2}
\\&~~~~~\quad
\quad \quad \times
\bigg( \int_{B_{2r}} |\nabla u| \d x  + r^{d-1} \max_{B_{2r}\setminus B_r} u \bigg) 
\\&\quad
\leq
\int_{B_{2r}\setminus B_r} u^{n} |\nabla \Delta u|^2 \d x
+C r^{d-6} \max_{B_{2r}\setminus B_r} u^{n+2}
\\&~~~\quad
+C r^{-d-4} \max_{B_{2r}} u^{n} \bigg(\max_{B_{2r}} u^{(4-n)/6} \int_{B_{2r}} |\nabla u^{(n+2)/6}| \d x\bigg)^2
\\&\quad
\leq
\int_{B_{2r}\setminus B_r} u^{n} |\nabla \Delta u|^2 \d x
+C r^{d-6} \max_{B_{2r}\setminus B_r} u^{n+2}
\\&~~~\quad
+C r^{2(d-6)/3} \max_{B_{2r}} u^{2(n+2)/3} \bigg(\int_{B_{2r}\setminus B_{\delta r}} |\nabla u^{(n+2)/6}|^6 \d x\bigg)^{1/3}
\\&~~~\quad
+C \delta^{5d/3} r^{2(d-6)/3} \max_{B_{2r}} u^{2(n+2)/3} \bigg(\int_{B_{\delta r}} |\nabla u^{(n+2)/6}|^6 \d x\bigg)^{1/3}.
\end{align*}
Using Lemma~\ref{BoundMaxUDegenerate} and Young's inequality, this yields the statement of the lemma.
\end{proof}

\section{Proof of Theorem \ref{MainResult}: space-time H\"older continuity for the 2D thin-film equation}

We first proof Theorem \ref{MainResult} in the case of stronger regularity of the initial datum stated in point b), i.e., $\nabla u_0 \in L^p(\mathbb{R}^d)$, for some $p>2$. We then perform the (short) adaptation needed to treat the baseline case $u_0 \in H^1(\mathbb{R}^d)$.
\begin{proof}[Proof of Theorem \ref{MainResult}, case b): $\nabla u_0 \in L^p(\mathbb{R}^d)$ for some $p>2$]

We proceed in three steps.

\emph{Step 1: obtaing power-law bounds $\int_{B_r}|\nabla u|^2 \d x \lesssim r^{\beta}$ using hole-filling estimate from Lemma \ref{LemmaHoldFill_TFE}}. 
Fix a time $t>0$, and a radius $r>0$. 
We iterate the hole filling estimate of Proposition \ref{LemmaHoldFill_TFE} (used in its notationally convenient form \eqref{hf_est_succinct_tfe}) over the hole-filling terms that progressively pop up on the right-hand-side. Informally, this reads as
\begin{align*}
& \tilt{t}{r}{\frac{1}{\delta} \cdot r} + \hfterms{0}{t}{r} \\
& \quad \leq \tilt{0}{\frac{2}{\delta} \cdot r}{\frac{2}{\delta} \cdot \frac{r}{2}} + (1-\theta) \cdot \hfterms{0}{t}{\frac{2}{\delta} \cdot r} \\
& \quad \leq \tilt{0}{\frac{2}{\delta} \cdot r}{\frac{1}{\delta}\cdot r} + (1-\theta) \cdot \left[ \tilt{0}{(\frac{2}{\delta})^2 \cdot r}{(\frac{2}{\delta})^2 \cdot \frac{r}{2}} + (1-\theta) \cdot \hfterms{0}{t}{(\frac{2}{\delta})^2 \cdot r} \right]
\\&\quad
\leq \dots .
\end{align*}
Set $\ratiodelta: = 2/\delta$. Iterating \eqref{hf_est_succinct_tfe} $K$ times, with $K := \lfloor \log_\Lambda (r^{-1}) \rfloor$, so as to relate the final  estimate to integrals over a ball of radius $R\sim 1$, we get
\begin{align}\label{iterated_hole_filling}
  \tilt{t}{r}{\frac{1}{\delta} \cdot r} \leq \sum_{k=1}^{K}{(1-\theta)^k\tilt{0}{\ratiodelta^k \cdot r}{\ratiodelta^k \cdot \frac{r}{2}}} + (1-\theta)^{K} \hfterms{0}{t}{\ratiodelta^{K}\cdot r}.
\end{align}
Since $\nabla u_0\in L^p$ with $p>2$, it is straightforward to deduce that $\tilt{0}{\ratiodelta^k \cdot r}{\ratiodelta^k \cdot \frac{r}{2}} \leq C (\ratiodelta^k r)^\gamma$, where $\gamma := 2(p-2)/p >0$. Additionally, $\hfterms{0}{t}{\ratiodelta^{K}\cdot r}$ is bounded thanks to Definition \ref{DefinitionEnergyDissipatingWeakSolution}. 
Plugging this in \eqref{iterated_hole_filling}, using Lemma \ref{LemmaDiffBsCs} grants
\begin{align}\label{r_polynom_est_tilt}
 \tilt{t}{r}{\frac{1}{\delta}\cdot r} = \frac{1}{2}\int_{B_r} |\nabla u(x,t)-b_{\delta^{-1} r}(t) \cdot x -c_{\delta^{-1} r}(t)|^2 \,\d x \leq C r^\beta,
\end{align}
where $\beta := \min\{-\log_\ratiodelta(1-\theta); \gamma\} > 0$. Using 
a telescopic sum argument, we obtain
\begin{align}\label{telescopic}
\left(\int_{B_r}{|b_{\delta^{-1} r} \cdot x + c_{\delta^{-1} r}|^2 \d x}\right)^{1/2} 
& \leq r (|b_{\delta^{-1} r}| r + |c_{\delta^{-1} r}|) \\
& \leq r\left( \sum_{k=0}^{K-1}|b_{\ratiodelta^k \delta^{-1} r}-b_{\ratiodelta^{k+1} \delta^{-1} r}||r| + |b_{\ratiodelta^{K} \delta^{-1} r}||r| \right) \nonumber \\
& \quad \quad + r\left( \sum_{k=0}^{K-1}|c_{\ratiodelta^k \delta^{-1} r}-c_{\ratiodelta^{k+1} \delta^{-1} r}| + |c_{\ratiodelta^{K} \delta^{-1} r}| \right). \nonumber
\end{align}
Plugging the estimates of Lemma \ref{LemmaDiffBsCs} in \eqref{telescopic}, and also using \eqref{r_polynom_est_tilt}, we deduce
\begin{align}\label{L2_bound_Bs_plus_Cs}
& \left(\int_{B_r}{|b_{\delta^{-1} r} \cdot x + c_{\delta^{-1} r}|^2 \d x}\right)^{1/2} \\
&
\nonumber
\quad \leq C r\left( \sum_{k=0}^{K-1}(\ratiodelta^k r)^{\beta/2-2}r + |b_{\ratiodelta^{K} \delta^{-1} r}||r| \right) + C r\left( \sum_{k=0}^{K-1}(\ratiodelta^k r)^{\beta/2-1} + |c_{\ratiodelta^{K}\delta^{-1} r}| \right)
\\&
\quad \leq C r^{\beta/2}\nonumber,
\end{align}
where, in the last inequality, we used the convergence of the sum $\sum_{k=0}^{\infty}{\ratiodelta^{k(\beta/2-1)}}$ (as $\beta/2-1<0$) and straightforward bounds on $b_{\ratiodelta^K \delta^{-1} r}$ and $c_{\ratiodelta^K \delta^{-1} r}$. The triangle inequality combined with \eqref{r_polynom_est_tilt} and \eqref{L2_bound_Bs_plus_Cs} gives the desired bound $\int_{B_r}|\nabla u|^2 \d x \leq C r^{\beta}$.

\emph{Step 2: Spatial H\"older continuity}.
The estimate $\int_{B_r}|\nabla u|^2 \d x \leq C r^{\beta}$ implies via the Poincar\'e inequality that $u$ belongs to the Campanato space $\mathcal{L}^{2+\beta,2}(\mathbb{R}^2)$, entailing H\"older continuity. For the reader's convenience, we briefly recall the classical argument:
By the Lebesgue differentiation theorem and (in the second step) the Poincar\'e inequality, we have for a.\,e.\ $x\in \mathbb{R}^2$
\begin{align*}
\bigg|u(x)-\dashint_{B_r(x)} u \,\d y\bigg| & \leq 
\sum_{k=0}^\infty\bigg|\dashint_{B_{2^{-k} r}(x)} u \,\d y
-\dashint_{B_{2^{-k-1}r}(x)} u \,\d y\bigg|
\\&
\leq C
\sum_{k=0}^\infty (2^{-k}r)^{2-d} \bigg(\int_{B_{2^{-k} r}(x)} |\nabla u|^2 \,\d y\bigg)^{1/2},
\end{align*}
which yields, in view of $\int_{B_{2^{-k} r}(x)} |\nabla u|^2 \,\d y\leq C (2^{-k}r)^{\beta}$ in the case $d=2$, that
\begin{align*}
&\bigg|u(x)-\dashint_{B_r(x)} u \,\d y\bigg|
\leq C r^{\beta/2}.
\end{align*}
For two points $x_1,x_2$, we obtain by setting $r:=|x_1-x_2|$
\begin{align*}
&|u(x_1)-u(x_2)|
\\&
\leq \bigg|u(x_1)-\dashint_{B_r(x_1)} u \,\d y\bigg|+\bigg|u(x_2)-\dashint_{B_r(x_2)} u \,\d y\bigg|
+\bigg|\dashint_{B_r(x_1)} u \,\d y-\dashint_{B_r(x_2)} u \,\d y\bigg|
\\&
\leq C r^{\beta/2} + C r^{\beta/2} + C r^{1-d/2} \bigg(\int_{B_{2r}(x_1)} |\nabla u|^2 \,\d y\bigg)^{1/2}
\\&
\leq C r^{\beta/2} = C|x_1-x_2|^{\beta/2},
\end{align*}
giving spatial H\"older continuity with parameter $\sigma_x = \frac{\beta}{2} = \min\left\{-\frac{\log_{\frac{2}{\delta}}(1-\theta)}{2}; \frac{p-2}{p}\right\}$.

\emph{Step 3: H\"older continuity in time}.
Consider a smooth test function $\eta$ supported in $B_r(x)$, with $\int \eta \,\d x =1$, and $|\nabla \eta| \leq C r^{-d-1}$. Then
\begin{align}\label{time_holder_est_prep}
&|u(x,t_1)-u(x,t_2)|
\\& \nonumber
\leq \bigg|u(x,t_1)-\int_{B_r(x)} \eta(y) u(y,t_1) \,\d y\bigg|
+\bigg|u(x,t_2)-\int_{B_r(x)} \eta(y) u(y,t_2) \,\d y\bigg|
\\&~~~ \nonumber
+\bigg|\int_{t_1}^{t_2} \int_{B_r(x)} \eta(y) \partial_t u(y,t) \,\d y \,\d t\bigg|
\\& \nonumber
\leq C r^{\sigma_x} +\bigg| \int_{t_1}^{t_2} \int u^n \nabla\Delta u \cdot \nabla \eta \,\d y \,\d t \bigg|
\\& \nonumber
\leq C r^{\sigma_x} +\bigg( \int_{t_1}^{t_2} \int u^n |\nabla\Delta u|^2 \,\d y \,\d t \bigg)^{1/2} \bigg( \int_{t_1}^{t_2} \int u^n |\nabla \eta|^2 \,\d y \,\d t \bigg)^{1/2}
\\& \nonumber
\leq C r^{\sigma_x} + C |t_2-t_1|^{1/2} r^{-d-1}  \bigg( \int_{t_1}^{t_2} \int u^n |\nabla\Delta u|^2 \,\d y \,\d t \bigg)^{1/2} \bigg( \dashint_{t_1}^{t_2} \int u^n \,\d y \,\d t \bigg)^{1/2}.
\end{align}
From Definition  \ref{DefinitionEnergyDissipatingWeakSolution}, we know that $\chi_{\{u > 0\}}u^{\frac{n}{2}}\nabla\Delta u \in L^2(\R^d \times [0,T))$, and this implies boundedness of the term in the first round bracket of the right-hand-side of \eqref{time_holder_est_prep}. The term in the second round bracket of \eqref{time_holder_est_prep} can be bounded using the fact that $L^\infty([0,T); H^1(\R^d))$ (again from Definition \ref{DefinitionEnergyDissipatingWeakSolution}) and the Sobolev embedding theorem $H^{1}(\mathbb{R}^d)\subset L^{n}(\mathbb{R}^d)$ for $n<3$.
This gives
\begin{align*}
|u(x,t_1)-u(x,t_2)| \leq C r^{\sigma_x} + C |t_2-t_1|^{1/2}r^{-d-1}.
\end{align*}
Optimizing in $r$ (i.e., taking $r\propto |t_2-t_1|^{1/(2(\sigma_x+d+1))}$) finally yields
\begin{align*}
|u(x,t_1)-u(x,t_2)|
\leq C |t_2-t_1|^{\frac{\sigma_x}{2(\sigma_x+d+1)}},
\end{align*}
entailing time H\"older continuity with parameter $\sigma_t := \frac{\sigma_x}{2(\sigma_x+d+1)}$. 
\end{proof}
\begin{proof}[Proof of Theorem \ref{MainResult}, case a): $u_0 \in H^1(\mathbb{R}^d)$] In this case, we can not pivot the hole-filling estimate on the initial condition, as it lacks sufficient regularity. However, Definition \ref{DefinitionEnergyDissipatingWeakSolution} tells us that $\int_{0}^{T}\int_{\mathbb{R}}|\nabla u^{\frac{n+2}{6}}|^6 \d x\d t \leq C \|u_0\|_{H^1(\mathbb{R}^d)}$ due to the Bernis-Gr\"un inequalities. This means that, for any $t\in (0,T)$, we can find $\tilde{t}\in (0,t/2)$ such that 
$
\|\nabla u^{\frac{n+2}{6}}(\cdot,\tilde{t})\|^6_{L^6(\mathbb{R}^2)}
\leq 2C \|u_0\|_{H^1(\mathbb{R}^d)} t^{-1}
$.
Using the H\"older inequality, the Sobolev embedding $H^1(\mathbb{R}^{d})\subset L^{4-n}(\mathbb{R}^{d})$, the energy dissipation, and the conservation of mass, we deduce
\begin{align}\label{regular_time}
\int_{\mathbb{R}^d}{|\nabla u(x,\tilde{t})|^3 \d x} & = \int_{\mathbb{R}^d}{u^{\frac{4-n}{2}}(x,\tilde{t})u^{\frac{n-4}{2}}(x,\tilde{t})|\nabla u(x,\tilde{t})|^3 \d x}
\\
& \leq C\|u(\cdot,\tilde{t})\|^{\frac{4-n}{2}}_{H^1(\mathbb{R}^d)} \|\nabla u^{\frac{n+2}{6}}(\cdot,\tilde{t})\|^3_{L^6(\mathbb{R}^2)}\nonumber\\
& \leq C(\|u_0\|_{H^1(\mathbb{R}^d)},n) t^{-1/2}\nonumber
\end{align}
and this means that $\nabla u(\cdot,\tilde{t})\in L^p(\mathbb{R}^d)$ with $p=3>2$. Therefore, we can follow the proof of spatial H\"older continuity in case b) simply by replacing $\{t_1;t_2\} = \{0;t\}$ with $\{t_1;t_2\} = \{\tilde{t};t\}$: of course, the H\"older continuity will only be local due the diverging term $t^{-1/2}$ in \eqref{regular_time}. The proof of time H\"older continuity (again, of local type due to the time singularity $t^{-1/2}$) is analogous to the one in case b).
\end{proof}

\appendix

\section{Relevant inequalities}

\begin{proposition}[The Bernis-Gr\"un inequality \cite{GruenBernis}, combined with the approximation argument in {\cite[Proof~of~Lemma~A.4]{DeNittiFischer}}]\label{BernisGruenIneq}
Let $d\in \{1,2,3\}$ and let $n\in (2-\sqrt{8/(8+d)},3)$. Let $\eta$ be a nonnegative smooth compactly supported weight.
Let $u\in L^1(\mathbb{R}^d)$ be a nonnegative function with $\nabla u^{(n+2)/6} \in L^6(\mathbb{R}^d)$ and $u^{n/2} \nabla \Delta u\in L^2(\mathbb{R}^d)$. Then we have $u\in H^3_{loc}(\{x \in \mathbb{R}^d:u(x)>0\})$ and there exists $C(d,n)$ such that
\begin{align}
\label{BernisGruenWeighted}
&\int \Big[|\nabla u^{(n+2)/6}|^6 + u^{n-2} |D^2 u|^2 |\nabla u|^2 + |D^3 u^{(n+2)/2}|^2 \Big] \eta^6 \,\d x
\\&
\nonumber
\leq C(d,n) \left[\int u^n |\nabla \Delta u|^2 \eta^6 \,\d x + \int u^{n+2} (|\nabla \eta|^6 + \eta^3 |D^2\eta|^3 + \eta^4 |D^3 \eta|^2) \,\d x\right].
\end{align}
\end{proposition}
Note that in Proposition~\ref{BernisGruenIneq} the second term on the left-hand side is to be understood as $u^{n-2} |D^2 u|^2 |\nabla u|^2=\tfrac{6^2}{(n+2)^2} u^{2(n-1)/3} |D^2 u|^2 |\nabla u^{(n+2)/6}|^2$.

\begin{lemma}[Morrey's inequality]\label{lemma_morrey_ineq}
Let $d\in \mathbb{N}$ and let $p$ satisfy $d < p \leq \infty$. Then any function $u \in W^{1,p}(B_1)$ satisfies (after possible redefinition on a set of vanishing Lebesgue measure)  $u \in C^{0,1-\frac{d}{p}}(B_1)$, and the estimate
\begin{align*}
\|u\|_{C^{0,1-\frac{d}{p}}(B_1)} \leq C \|u\|_{W^{1,p}(B_1)}
\end{align*}
holds for some $C=C(d,p)$.
\end{lemma}

\section{Miscellaneous}\label{appendix_Bs_Cs}

\begin{lemma}\label{finite_2_3_good_times}
Let $u$ be a solution to the thin-film equation \eqref{tfe} in the sense of Definition~\ref{DefinitionEnergyDissipatingWeakSolution}.
Let $t$ be a good time for radius $r$ in the sense of Definition~\ref{DefnGoodBad}. Then $\int_{B_r}{|D^2u|^2}\d x$ and $\int_{B_r}{|D^3u|^2}\d x$ are finite. 
\end{lemma}
\begin{proof}
If $\inf_{x\in B_r}u =0$, then $u\equiv 0$ on $B_r$ by the good time condition, thus the claim is trivial. If instead $\inf_{x\in B_r}u >0$, we can write
\begin{align}
\int_{B_r}{|D^2u|^2}\d x & \leq (\inf_{x\in B_r}u)^{-(n+\alpha-1)}\int_{B_r}{u^{n+\alpha-1}|D^2u|^2\d x},
\label{bound_full_second_der}
\\
\int_{B_r}{|D^3u|^2}\d x & \leq (\inf_{x\in B_r}u)^{-n}\int_{B_r}{u^n|D^3u|^2\d x}. 
\label{bound_full_third_der}
\end{align}
The right-hand-side of \eqref{bound_full_second_der} is finite thanks to Definition \ref{DefinitionEnergyDissipatingWeakSolution}, point b). The right-hand-side of \eqref{bound_full_third_der} is finite due to \eqref{control_un_full_third_derivative}, Definition \ref{DefinitionEnergyDissipatingWeakSolution} (point a)), and the Sobolev embedding $H^1(\mathbb{R}^d)\subset L^{n+2}(\mathbb{R}^d)$, $n<3$.
\end{proof}

\begin{lemma}\label{lemma_poincare}
Let $r>0$. Let $u\in H^3(B_{2r})$ and let $b_r$ and $c_r$ be as defined in \eqref{defb}--\eqref{defc}. Let $t$ be a good time for radius $2r$. Then we have
    \begin{align}
    \int_{B_{2r}\setminus B_{r}}{|\nabla u - b_r \cdot x -c_r|^2}
& \leq C r^2 \int_{B_{2r}\setminus B_r}|D^2 u - b_r|^2 \d  x, \label{poincare_nabla}\\
\int_{B_{2r}\setminus B_r}|D^2 u - b_r|^2\d x & \leq C r^2 \int_{B_{2r}\setminus B_r}|D^3 u|^2 \d x. \label{poincare_d2}
    \end{align}
The same estimates also hold true if the annulus $B_{2r}\setminus B_r$ is replaced in all integrals by the ball $B_{2r}$. Furthermore, for $d=2$ we also have the Poincar\'e-Sobolev inequality
\begin{align}
\label{PoincareSobolev}
\int_{B_{2r}}|D^2 u - b_r|^2\d x & \leq C \bigg(\int_{B_{2r}}|D^3 u| \d x\bigg)^2.
\end{align}
\end{lemma}
\begin{proof}
We only prove \eqref{poincare_nabla}-\eqref{poincare_d2} on the annulus $B_{2r}\setminus B_r$; the proof of \eqref{poincare_nabla}-\eqref{poincare_d2} on the ball $B_{2r}$ is analogous, as is the proof of the Poincar\'e-Sobolev inequality \eqref{PoincareSobolev}.

Note that, if $t$ is a good time for $2r$, all integrals in \eqref{poincare_d2}--\eqref{poincare_nabla} are well defined thanks to Lemma \ref{finite_2_3_good_times}. To prove \eqref{poincare_d2}, set $\overline{b}_r := \dashint_{B_{2r}\setminus B_r}{D^2 u \d x}$, i.e., let $\overline{b}_r$ be the average of $D^2 u$ on $B_{2r}\setminus B_r$. Abbreviate 
$$
m(x):= \frac{\tilde{\eta}(\frac{x}{r})}{\int_{B_{2r}\setminus B_r}{\tilde{\eta}(\frac{x}{r})}} - \frac{1}{|B_{2r}\setminus B_r|}.
$$
The triangle inequality and integration by parts in the definition of $b_r$ grant   \begin{align*}
& \int_{B_{2r}\setminus B_r}|D^2 u - b_r|^2\d x \\
& \leq 2\int_{B_{2r}\setminus B_r}|D^2 u - \overline{b}_r|^2\d x +2 \int_{B_{2r}\setminus B_r}|b_r - \overline{b}_r|^2\d x \\
& \leq 2\int_{B_{2r}\setminus B_r}|D^2 u - \overline{b}_r|^2\d x + 2|B_{2r}\setminus B_r||b_r - \overline{b}_r|^2 \\
& \leq C \int_{B_{2r}\setminus B_r}|D^2 u - \overline{b}_r|^2\d x + C r^2\left|\int_{B_{2r}\setminus B_r}{D^2 u : m(x) \d x}\right|^2 \\
& \leq C \int_{B_{2r}\setminus B_r}|D^2 u - \overline{b}_r|^2\d x + C r^2\left|\int_{B_{2r}\setminus B_r}{(D^2 u - \overline{b}_r ) : m(x) \d x}\right|^2, 
\end{align*}
where we have used the fact that $m(x)$ has mean zero. By $|m(x)|\leq C r^{-2}$, we thus obtain the bound $\int_{B_{2r}\setminus B_r}|D^2 u - b_r|^2\d x \leq C \int_{B_{2r}\setminus B_r}|D^2 u - \overline{b}_r|^2\d x$. Applying the Poincar\'e inequality on the right-hand-side of this estimate gives \eqref{poincare_d2}.
 
As for \eqref{poincare_nabla}, define $\overline{c}_r:= \dashint_{B_{2r}\setminus B_r}{(\nabla u - b_r \cdot x)\d x}$. Then 
\begin{align*}
    & \int_{B_{2r}\setminus B_{r}}{|\nabla u - b_r \cdot x -c_r|^2} \d x \\
    & \quad \leq C \int_{B_{2r}\setminus B_{r}}{|\nabla u - b_r \cdot x -\overline{c}_r|^2} \d x + C \int_{B_{2r}\setminus B_{r}}{|c_r -\overline{c}_r|^2} \d x.
\end{align*}
Since $m(x)$ is radially symmetric and with mean zero, and $b_r\cdot x$ is radially anti-symmetric, we deduce $|\overline{c}_r - c_r| = |\int_{B_{2r}\setminus B_{r}}{m(x)(\nabla u(x)-b_r \cdot x - \overline{c}_r)\d x}|$. Inequality \eqref{poincare_nabla} follows promptly using Poincar\'e and \eqref{poincare_d2}.
\end{proof}

\begin{lemma}\label{LemmaDiffBsCs}
Let $d=2$.
Let $b_r$ and $c_r$ be as defined in \eqref{defb}--\eqref{defc}. 
Let $\delta \in (0,1]$. Then
\begin{align}
|b_r - b_{\frac{\delta r}{2}}|^2 & \leq C r^{-4}\int_{B_{2r}}{\frac{1}{2}|\nabla u - b_r\cdot x - c_r|^2} \d x ,\label{DiffBs_est}\\
|c_r - c_{\frac{\delta r}{2}}|^2 & \leq C r^{-2}\int_{B_{2r}}{\frac{1}{2}|\nabla u - b_r\cdot x - c_r|^2} \d x \label{DiffCs_est}.
\end{align}
\end{lemma}
\begin{proof}
Define 
\begin{align*}
m(x) := -\frac{1}{\int_{B_{2r}} \tilde{\eta}(\frac{x}{r})\d x} \tilde{\eta}\left(\frac{x}{r}\right) 
 +\frac{1}{\int_{B_{\delta r/2}} \tilde{\eta}(\frac{x}{\delta r/2})\d x}\tilde{\eta}\left(\frac{x}{\delta r/2}\right).
\end{align*}
Since $\tilde{\eta}(\frac{\cdot}{\delta r/2})$ and $\tilde \eta(\tfrac{\cdot}{r})$ are supported on $B_{2r}$, we have
\begin{align}\label{BsDiff}
& b_r - b_{\frac{\delta r}{2}} = -
 \int_{B_{2r}} \nabla m(x) \otimes \nabla u \d x.
\end{align} 
In \eqref{BsDiff} we can replace $\nabla u$ with $\nabla u - b_r \cdot x - c_r$: this is the case since (once integrating by parts by using the fact $m(x)$ vanishes at the boundary of $B_{2r}$), $c_r$ does not contribute anything as it gradient is trivially zero, and the constant matrix $b_r$ is integrated against the average-zero function $m(x)$. 
  Estimate \eqref{DiffBs_est} then follows using H\"older's inequality \eqref{BsDiff} and the bounds $|m|\leq C r^{-d}$, $|\nabla m|\leq C r^{-d-1}$. To obtain \eqref{DiffCs_est} we observe that $c_r - c_{\frac{\delta r}{2}} = \int_{B_{2r}} m(x) \otimes \nabla u \d x = \int_{B_{2r}} m(x) \otimes (\nabla u -b_r \cdot x - c_r )\d x $, where we could add $c_r$ as its integrated against the average-zero function $m(x)$, as well as $b_r \cdot x$ (since it is radially anti-symmetric, and it's integrated against the radially symmetric $m(x)$). Using the H\"older inequality and the estimates on $m$ grants \eqref{DiffCs_est}.  
\end{proof}


\bibliographystyle{plain}
\bibliography{thinfilm.bib}

\end{document}